\documentclass[10pt]{amsart}
\usepackage{amssymb,amsmath,amsthm,amscd}
\usepackage{epsfig, comment,enumerate, MnSymbol}
\usepackage[all]{xy}

\newtheorem{thm}{Theorem}
\newtheorem{cor}[thm]{Corollary}
\newtheorem{lem}[thm]{Lemma}

\newtheorem{theorem-question}[thm]{Theorem-Question}

\theoremstyle{definition}
\newtheorem{defn}[thm]{Definition}

\newenvironment{exafont}{\begin{bf}}{\end{bf}}

\newcommand{\bfb}{\mathfrak{b}}\newcommand{\bfc}{\mathbf{c}}\newcommand{\bft}{\mathbf{t}}
\newcommand{\bfy}{\mathbf{y}}\newcommand{\bfY}{\mathbf{Y}}

\newcommand{\bbO}{\mathbb{O}}\newcommand{\fO}{\mathfrak{O}}
\newcommand{\bbH}{\mathbb{H}}

\newcommand{\cE}{\mathcal{E}}\newcommand{\bbE}{\mathbb{E}}

\newcommand{\NN}{\mathbb{N}}
\newcommand{\bbT}{\mathbb{T}}

\DeclareMathOperator{\End}{End}
\DeclareMathOperator{\Hom}{Hom}

\DeclareMathOperator{\Ext}{Ext}

\DeclareMathOperator{\im}{Im}
\DeclareMathOperator{\Ker}{Ker}\DeclareMathOperator{\Coker}{Coker}

\DeclareMathOperator{\ml}{-mod}

\DeclareMathOperator{\gr}{-gr}

\DeclareMathOperator{\perfl}{-perf}
\DeclareMathOperator{\perfr}{perf-}

\DeclareMathOperator{\perf}{-perf}

\title{Koszul dual $2$-functors and extension algebras of simple modules for $GL_2$}

\author{Vanessa Miemietz and Will Turner}
\date{}

\begin{document}

\begin{abstract}
Let $p$ be a prime number.
We compute the Yoneda extension algebra of $GL_2$ over an
algebraically closed field of characteristic $p$ by developing a theory of Koszul duality for a certain class of $2$-functors,
one of which controls the category of rational representations of $GL_2$ over such a field.
\end{abstract}


\maketitle

\tableofcontents
\setlength{\parskip}{8pt}
\section{Introduction.}\label{intro}

Let $F$ be an algebraically closed field of characteristic $p>0$.
Let $G = GL_2(F)$ denote the group of $2 \times 2$ invertible
matrices over $F$. Let $\mathcal{L}$ denote a complete set of
irreducible objects in the category $G \ml$ of rational
representations of $G$. The object of this paper is to give an explicit description of the Yoneda extension algebra
$$\bfY = \bigoplus_{L,L' \in \mathcal{L}} \Ext_{G \ml}^\bullet(L,L')$$ of $G$.
The strongest previous results in this direction were obtained by A. Parker, who outlined an intricate algorithm to compute the dimension of $\Ext_{G\ml}^n(L,L')$ for $L, L' \in \mathcal{L}$ and $n \geq 0$ \cite{AP}.
Our approach is quite different. We develop a theory of homological duality for certain algebraic operators $\mathbb{O}$ which we
introduced previously in our study of the category of rational representations of $G$. Heuristically speaking, these operators (which are indeed $2$-functors on a certain $2$-category) form a way of "applying" one algebra to another to produce a bigger algebra (dragging certain bimodules around). While $G\ml$ is not Koszul itself, it is obtained as an iterated application of Koszul algebras to the ground field, and in this article we show that iteratively applying their Koszul duals (together with the appropriate bimodules) will produce a dg algebra  whose homology is the desired extension algebra.
We use this to give a combinatorial description of $\bfY$ as an algebra, including the multiplicative structure, described by explicit formulae how to multiply basis elements. While in the end, we do arrive at an explicit combinatorial description of $\bfY$, we would like to stress that the iterated dg algebra, as whose homology we obtain the extension algebra, keeps track of all the homotopical information, and hence in principle incorporates the (nontrivial) $A_\infty$-structure of the Yoneda extension algebra of $GL_2$. To investigate this $A_\infty$-structure would in itself be an interesting problem.

The category $G \ml$ has countably many blocks, all of which are
equivalent. Therefore, the algebra $\bfY$ is isomorphic to a direct sum
of countably many copies of $\bfy$, where $\bfy$ is the Yoneda extension
algebra of the principal block of $G$. Our problem is to compute
$\bfy$. In the remainder of this section, we will introduce the main objects in the explicit computations, so we can state our main theorem, Theorem \ref{Yoneda}. We will then deduce the description of a multiplicative (up to prescribed sign) basis from this, and finish with some remarks. In Section \ref{example}, we will give an explicit description of an example illustrating the monomial basis. In Section \ref{kellerkoszul}, we  introduce the $2$-category in which our theory is developed and extend notions of homological duality to this setting. In Section \ref{ops}, we introduce our $2$-functors, before studying their behaviour under the homological dualities from  Section \ref{kellerkoszul} in Section \ref{respect}. We then recall the results needed from the representation theory of $GL_2(F)$ in Section \ref{GL2}. Section \ref{reduction} then uses the developed theory to reduce the computation to the computation of the homology of an explicit tensor algebra. The latter is then computed in Section \ref{explicit}, which is largely independent of the preceding sections.

Let $\Pi$ denote the
preprojective algebra of bi-infinite type $A$: it is the path
algebra of the quiver
$$
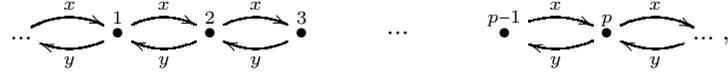
\xymatrix{ ... \ar@/^/[r]^{x}
& \ar@/^/[l]^{y} \overset{1}{\bullet} \ar@/^/[r]^{x}
&\overset{2}{\bullet}\ar@/^/[l]^{y} \ar@/^/[r]^{x} &\ar@/^/[l]^{y}
\overset{3}{\bullet} &\cdots &\overset{p-1}{\bullet}\ar@/^/[r]^{x} &
\ar@/^/[l]^{y} \overset{p}{\bullet} \ar@/^/[r]^{x} & \ar@/^/[l]^{y}
...},$$ modulo relations $xy-yx$, and  is naturally
$\mathbb{Z}_+$-graded so that paths are homogeneous with degree
given by their length. 

This algebra has a natural basis given by
$$\mathcal{B}_{\Pi} = \{ (s,\alpha,\beta) \quad | \quad s \in \mathbb{Z}, \alpha, \beta \in \mathbb{N}\}.$$
Here, basis elements of $\Pi$ correspond to paths in the quiver of $\Pi$; the element $(s,\alpha,\beta)$ corresponds to the element $x^\alpha y^\beta e_s$ in the path algebra. The target $t$ of such a path can be computed via the
formula $t-s = \alpha - \beta$. The product in $\Pi$ of two such basis elements is defined by the following formula:
$$ (s,\alpha,\beta)\cdot (s',\alpha',\beta') = \left\{ \begin{array}{ll}
(s, \alpha+\alpha',\beta+\beta') & \mbox{if $s' = s+\alpha-\beta$}; \\
0 & \mbox{otherwise}.
\end{array} \right.$$
In the following, we will call such a basis for an algebra monomial, meaning the basis contains a complete set of primitive
idempotents and the product of any two of its elements is either equal to ($\pm$) another basis element or equal to
zero.

For $b = (s,\alpha,\beta) \in \mathcal{B}_{\Pi}$, we define the
degree of $b$ to be $|b| = \alpha+\beta$.

Let $\Pi_{\leq p}$ denote the subalgebra of $\Pi$ generated by
arrows which begin and end at vertices indexed by $i \leq p$. Let
$\Omega$ denote the quotient of $\Pi_{\leq p}$ by the ideal generated by
idempotents $e_i$ corresponding to vertices $i \leq {0}$, that is $\Omega$ is given by quiver and relations as $$\Omega = 
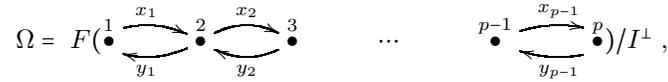
\xymatrix{F(\overset{1}{\bullet}
\ar@/^/[r]^{x_1} &\overset{2}{\bullet}\ar@/^/[l]^{y_1}
\ar@/^/[r]^{x_2} &\ar@/^/[l]^{y_2} \overset{3}{\bullet} &\cdots
&\overset{p-1}{\bullet}\ar@/^/[r]^{x_{p-1}} & \ar@/^/[l]^{y_{p-1}}
\overset{p}{\bullet})/I^\perp },$$ where $I^\perp=(x_l y_l - y_{l+1}
x_{l+1} , y_1x_1 \mid 1 \leq l \leq p-2)$.. 

The Loewy structure of the $p$th projective $\Omega e_p$ of this algebra is given by
{\footnotesize 
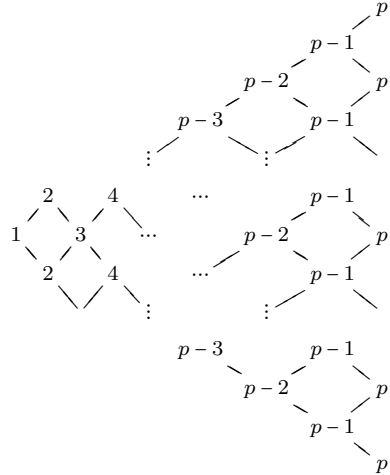
$$\xymatrix@=3pt{
&&&&&&&& p \ar@{-}[dl] \\
&&&&&&& p-1\ar@{-}[dl] \ar@{-}[dr]  &\\
&&&&&& p-2 \ar@{-}[dl] \ar@{-}[dr] && p\ar@{-}[dl] \\
&&&&& p-3 \ar@{-}[dl] \ar@{-}[dr] && p-1 \ar@{-}[dl] \ar@{-}[dr] &\\
&&&&\vdots&&\vdots&&\\
& 2 \ar@{-}[dl] \ar@{-}[dr]&& 4 \ar@{-}[dl] \ar@{-}[dr]&& \cdots && p-1 \ar@{-}[dl] \ar@{-}[dr]& \\
1 \ar@{-}[dr]&& 3 \ar@{-}[dl] \ar@{-}[dr]&& \cdots &&p-2 \ar@{-}[dl] \ar@{-}[dr]&& p\ar@{-}[dl] \\
& 2 \ar@{-}[dr]&& 4 \ar@{-}[dl] \ar@{-}[dr]&& \cdots && p-1\ar@{-}[dl] \ar@{-}[dr] & \\
&&&&\vdots&&\vdots&&\\
&&&&& p-3 \ar@{-}[dr]&& p-1 \ar@{-}[dl] \ar@{-}[dr]&\\
&&&&&& p-2 \ar@{-}[dr]&& p\ar@{-}[dl] \\
&&&&&&& p-1 \ar@{-}[dr]&\\
&&&&&&&& p
}$$}
and the remaining projectives are the corresponding submodules of this module.

A monomial basis for $\Omega$ is given by the subset
$$\mathcal{B}_{\Omega} = \{ (s,\alpha,\beta) \in \mathcal{B}_{\Pi} \quad | \quad 1 \leq s \leq p, \alpha \leq \beta+p-s, \beta \leq s-1\}$$
of $\mathcal{B}_{\Pi}$. When multiplying elements of
$\mathcal{B}_{\Omega}$, we use the multiplication rule for
$\mathcal{B}_{\Pi}$, with the caveat that products of basis elements
are zero if their product in $\mathcal{B}_{\Pi}$ does not belong to
$\mathcal{B}_{\Omega}$.

Let $\Theta$ denote the quotient of $\Pi$ by the ideal generated by idempotents corresponding to
vertices $i \leq 0$ or $i \geq p$. Alternatively, we can describe $\Theta$ as the quotient of $\Omega$ by the ideal
generated by vertex $p$, i.e.\ $\Theta = \Omega/\Omega e_p\Omega$. We will denote the natural quotient map $\Omega \to \Theta$ by $\pi$.
If $p=5$, the Loewy structure is given by
{\footnotesize 
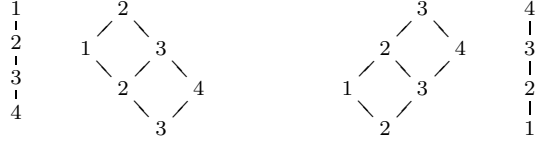
$$\xymatrix@=3pt{1\ar@{-}[d]\\2\ar@{-}[d]\\3\ar@{-}[d]\\ 4
}  \qquad
\xymatrix@=5pt{
& 2 \ar@{-}[dl] \ar@{-}[dr]&&&&\\
1 \ar@{-}[dr]&& 3\ar@{-}[dl] \ar@{-}[dr] &&&\\
& 2 \ar@{-}[dr]&& 4 \ar@{-}[dl] &&\\
&&3&&&&
} \qquad 
\xymatrix@=5pt{&& 3 \ar@{-}[dl] \ar@{-}[dr]&\\
& 2 \ar@{-}[dl] \ar@{-}[dr]&& 4 \ar@{-}[dl]\\
1\ar@{-}[dr]&& 3 \ar@{-}[dl]&\\
& 2 &&
}\qquad 
\xymatrix@=5pt{4\ar@{-}[d]\\3\ar@{-}[d]\\2\ar@{-}[d]\\1
}
$$}
which conveys the general picture.

A monomial basis for $\Theta$ is given by the
subset
$$\mathcal{B}_{\Theta} = \{(s,\alpha,\beta) \in \mathcal{B}_{\Pi} \quad | \quad 1 \leq s \leq p-1, \alpha \leq p-s-1, \beta \leq s-1
\}$$
of $\mathcal{B}_{\Pi}$.

We denote by $\sigma$ the involution of $\Pi$ which sends vertex $i$ to vertex $p-i$ and exchanges $x$ and $y$. This induces an action of $\sigma$ on $\Theta$ which, on a basis element $(s,\alpha,\beta) \in \mathcal{B}_{\Theta}$, is given by the formula
$$\sigma(s,\alpha,\beta) = (p-s,\beta,\alpha).$$ 
Let $\Theta^{\sigma}$ denote the $\Omega$-$\Omega$ bimodule given by lifting along $\pi$ the $\Theta$-$\Theta$-bimodule obtained by twisting the regular $\Theta$-bimodule on the right by $\sigma$.

We now define
$$\Lambda = \bbT_{\Omega}(\Theta^{\sigma}) \otimes F[\zeta],$$
to be the tensor product over $F$ of the tensor algebra over
$\Omega$ of $\Theta$ twisted on the right by $\sigma$, with a
polynomial algebra in a single variable. Note that as a vector space $(\Theta^\sigma)^{\otimes n}$ is just isomorphic to $\Theta$.

A monomial basis $\mathcal{B}_\Lambda$ for $\Lambda$ is given by
$$\mathcal{B}_\Lambda = \{(b,n,h) | n, h \in
\mathbb{N}, b \in \mathcal{B}_{\Omega} \textrm{ if } n=0, b
\in \mathcal{B}_{\Theta} \textrm{ if } n>0\}.$$
The product is given by
$$(b,n,h)(b',n',h') = \left\{\begin{array}{ll}
(bb',n+n',h+h') & \mbox{if $n$ is even}; \\
(b\sigma(b'),n+n',h+h') & \mbox{ if $n$ is odd}.
\end{array} \right.$$
Here, a basis element $(b,n,h)$  belongs to the component $(\Theta^{\sigma})^{\otimes n} \otimes \zeta^h$ of $\Lambda$.

We define a trigrading
$\Lambda = \bigoplus_{i,j} \Lambda^{ijk}$ on $\Lambda$ as follows:
the $i$-grading is defined by placing $\Omega$ in $i$-degree $0$ and
$\Theta^{\sigma}$ and $\zeta$ in $i$-degree $1$; the $j$-grading is defined
by grading elements of $\Omega$ and $\Theta^{\sigma}$ according to
path length and placing $\zeta$ in $j$-degree $p$; the $k$-grading is defined
by grading elements of $\Omega$ and $\Theta^{\sigma}$ according to
path length and placing $\zeta$ in $k$-degree $p-1$.
For a basis element $\beta=(b,n,h)$ this yields the $ijk$-degree $$(|\beta|_i,|\beta|_j,|\beta|_k)=(n+h, ph + |b|, |b|+(p-1)h).$$

Now suppose $\Gamma = \bigoplus_{i,j,k \in \mathbb{Z}} \Gamma^{ijk}$ is a $\mathbb{Z}$-trigraded algebra. We have a combinatorial operator $\fO_\Gamma$
which acts on the collection of bigraded algebras $\Delta$ after the formula
$$\fO_\Gamma(\Delta)^{ik} = \bigoplus_{j, k_1+k_2=k} \Gamma^{ijk_1} \otimes \Delta^{jk_2},$$
where $\otimes$ denotes the super tensor product i.e.\ multiplication is given by
$$(\gamma_1\otimes \delta_1)(\gamma_2\otimes \delta_2) = (-1)^{|\gamma_2|_k |\delta_1|_k} (\gamma_1\gamma_2\otimes \delta_1\delta_2)$$
where $\gamma_1,\gamma_2 \in \Gamma, \delta_1,\delta_2 \in \Delta$ and $|\cdot|_k$ denotes the $k$-degree of the corresponding element.

We consider the field $F$ as a trigraded algebra concentrated in degree $(0,0,0)$. We have a natural embedding of bigraded algebras $F
\rightarrow \Lambda$, which sends $1$ to the idempotent $e_1\otimes 1$, where $e_1 \in \Omega$ is the idempotent corresponding to vertex $1$. This embedding lifts to a morphism of operators $\fO_F \rightarrow \fO_\Lambda$. We have $\fO_F^2 = \fO_F$. Putting these together, we obtain a sequence of operators
$$\fO_F \rightarrow \fO_F \fO_\Lambda \rightarrow \fO_F \fO_\Lambda^2 \rightarrow ...$$
which, applied to the bigraded algebra $F[z]$ with $z$ placed in $jk$-degree $(1,0)$, gives a sequence of algebra embeddings
$$\lambda_1 \rightarrow \lambda_2 \rightarrow \lambda_3 \rightarrow ...,$$
where $\lambda_q = \fO_F \fO_\Lambda^q(F[z])$. Taking the union of the algebras in this sequence gives us an algebra $\lambda$. Our
main theorem is the following:

\begin{thm} \label{Yoneda}
The algebra $\bfy$ is isomorphic to $\lambda$ as a bigraded algebra.
\end{thm}

The proof of this theorem will occupy the rest of the article. For now, we would like to use the monomial basis given for $\Lambda$ to describe a monomial basis for $\lambda$. We will do this by describing a monomial basis for each $\lambda_q$.

Firstly, we put

$$\mathcal{B}_q = \{[\beta_1,...,\beta_q,l] | \beta_i \in \mathcal{B}_\Lambda, l \in \mathbb{N}\},$$

which, identifying the tuple $[\beta_1,...,\beta_q,l]$ with $\beta_1\otimes \cdots \otimes \beta_q \otimes z^l$ forms a monomial basis of the super tensor product 
$\Lambda^{\otimes q} \otimes F[z]$, where again multiplication is given by
$$[\beta_1,...,\beta_q,l][\beta_1',...,\beta_q',l'] =
(-1)^{\sum_{s'<s} |\beta_s|_k|\beta'_{s'}|_k } [\beta_1\beta_1',...,\beta_q\beta_q',l+l'].$$

By construction, our algebra $\lambda_q = \fO_F \fO_\Lambda^q(F[z])$ is a subalgebra of $\Lambda^{\otimes q} \otimes F[z]$ consisting of linear combinations of those tuples where the $i$-degree of $\beta_1$ is zero, the $i$-degree of $\beta_s$ equals the $j$-degree of $\beta_{s-1}$ for $s=2,\dots, q$, and $l = |\beta_q|_j$.
Hence defining a height function on $\mathcal{B}_q$ via 
$$[\beta_1,...,\beta_q,l] \mapsto (|\beta_{1}|_i , |\beta_{2}|_i-|\beta_{1}|_j,  |\beta_{3}|_i-|\beta_{2}|_j, \dots,  |\beta_{q}|_i-|\beta_{q-1}|_j, l-|\beta_{q}|_j)$$
we obtain a basis $\mathcal{B}_{\lambda_q}$ of $\lambda_q$ as the subset of $\mathcal{B}_q$ of height $(0,\dots,0)$.

Letting $\hat{1} = ((1,0,0),0,0)$ denote the element of ijk-degree $(0,0,0)$ in $\mathcal{B}_\Lambda$ of  which corresponds to the idempotent in $\Omega$ indexed by the vertex $1$, we obtain an embedding of $\mathcal{B}_q$ in $\mathcal{B}_{q+1}$ which
takes $(\beta_1,...,\beta_q,l)$ to $(\hat{1},\beta_1,...,\beta_q,l)$. This embedding is multiplicative
and height preserving, and restricts to the embedding of $\lambda_q$ into $\lambda_{q+1}$ described earlier.

Hence, letting $\mathcal{B}$ denote the union of the sequence of embeddings
$$\mathcal{B}_1 \rightarrow \mathcal{B}_2 \rightarrow \mathcal{B}_3 \rightarrow ...$$
and $\mathcal{B}_\lambda$ the subset of $\mathcal{B}$ consisting of elements of height $0$, we arrive at the following:

\begin{thm} The basis
$\mathcal{B}_\lambda$, with product obtained by restriction from
$\mathcal{B}$, forms a monomial basis for $\lambda$.
\end{thm}

We conclude this section with some remarks.

{\bf The Yoneda grading.} While our algebras are multigraded, the grading of most interest for representation theorists is the grading by extension degree, or Yoneda grading. An element of $\mathcal{B}$ represents an extension of degree $d$ if its total $k$-degree is $d$. Note that the total $k$-degree of an element $[\beta_1,...,\beta_q,l] $ for $\beta_i = (b_i,n_i,h_i)$ is given by $$d=\sum_{i=1}^ q |b_i| + (p-1)h_i.$$ Since for $1 \leq i \leq q-1$ we have $|b_i|+ph_i = n_{i+1}+h_{i+1}$,  the sum rearranges to
$$d= |b_q|+h_q + \sum_{i=1}^q n_i = l + \sum_{i=1}^q n_i$$
where we have used that $n_1 = 0$.

{\bf Polytopes.} The basis $\mathcal{B}_\lambda$ is infinite.
However, it is a union of monomial bases $\mathcal{B}_{\lambda_q}$
for $\lambda_q$, whose elements are in natural one-one
correspondence with elements of finite lattice polytopes $P_q$ of
dimension $4q$.
We can build this up as follows: by definition elements of $\mathcal{B}_\Omega$ and
$\mathcal{B}_\Theta$ are indexed by lattice elements of a finite
polytope in $\mathbb{Z}^3$; elements of $\mathcal{B}_\Lambda$ are
indexed by lattice elements of a polytope in $\mathbb{Z}^5$;
elements of $\mathcal{B}_q$ are indexed by lattice elements of an
infinite nonconvex polytope in $\mathbb{Z}^{5q+1}$; elements of
$\mathcal{B}_q$ of height $0$ are indexed by elements of a polytope
in $\mathbb{Z}^{5q+1 - (q+1)}$ since they are the intersection of
$\mathcal{B}_q$ with the kernel of a linear surjection from
$\mathbb{Z}^{5q+1}$ to $\mathbb{Z}^{q+1}$; elements of a monomial
basis $\mathcal{B}_{\lambda_q}$ for $\lambda_q$ are therefore
indexed by lattice elements of a polytope $P_q$ in
$\mathbb{Z}^{4q}$; the polytope $P_q$ is finite because $\lambda_q$
is finite dimensional.
Computing dimensions Ext-groups now amounts to counting points in a polytope with a certain boundary condition (the total $k$-degree). 

{\bf Recipe for application to $GL_2$.}
Suppose you want to know the dimension of $\Ext^d(L(\nu), L(\mu))$, where $\nu, \mu$ are two highest weights of simple modules for $GL_2(F)$, each given by a non-negative integer. First you need to determine whether they are in the same block.  Weights in a block are linearly ordered, so we have a natural order-preserving bijection between weights in a given block and the natural numbers $1,2,\dots$. The natural number associated to a weight $\nu$ will be denoted by $m_\nu$. 
The algorithm to determine whether two weights are in the same block is sketched in \cite[Section 1]{AP}, and we repeat and refine it here for the reader's convenience. Write $\nu = \sum_{i=1}^r a_ip^i$, $\mu = \sum_{i=1}^r b_ip^i$ with $0 \leq a_i,b_i \leq p-1$, and $r$ the greater one of the $p$-adic valuations of $\nu$ and $\mu$.
If at least one of $a_0,b_0$ is not equal to $p-1$, the weights $\nu, \mu$ are in the same block if and only if either $\sum_{i=1}^r (a_i-b_i)p^{i-1}\equiv 0 \mod 2$ and $a_0=b_0$ or $\sum_{i=1}^r (a_i-b_i)p^{i-1}\equiv 1 \mod 2$ and $a_0=p-2-b_0$. If $\nu, \mu$ are in the same block, the natural numbers to which they correspond in the linear order are given by $m_\nu:=1+\sum_{i=1}^r (a_i)p^{i-1}$ and $m_\mu:=1+\sum_{i=1}^r (b_i)p^{i-1}$ respectively.
If both $a_0, b_0$ equal $p-1$, take the first index (say $s$) such that one of $a_s,b_s$ differs from $p-1$, and repeat the above algorithm with $\bar\nu = 
\sum_{i=s}^r a_ip^{i-s}$, $\bar\mu = \sum_{i=1}^r b_ip^{i-s}$.

If $\nu$ and  $\mu$ are in the same block corresponding to natural numbers $m_\nu$ and $m_\mu$ respectively, we then find the number of occurrences of the $m_\nu$th simple in the $m_\mu$th projective for $\lambda_r$, which have $k$-degree $d$.

For example, simple representations in  the principal bock have highest weights  of the form $\nu_n = 2(n-1)p$ or $\mu_n = (2n-1)p +p-2$ for all $n \geq 1$. These get assigned values $m_{\nu_n }=2n-1$ and $m_{\mu_n }=2n$ respectively in our notation. Then we write these numbers as $a_1+\sum_{i=2}^{q} (a_i-1)p^{i-1}$ with $1 \leq a_i \leq p$ and the idempotent at the corresponding vertex (i.e. the identity morphism on the correponding simple)  is given by the basis element $[((a_q,0,0),0,0), \dots , ((a_1,0,0),0,0),0] \in \mathcal{B}_{\Lambda_q}$. Any element of $\Ext^1$ will be given by a basis element of the form 
\begin{equation*}\begin{split}&[((a_q,0,0),0,0), \dots ,((a_2,0,0),0,0), ((a_1,1,0),0,0),1] ,\\&
[((a_q,0,0),0,0), \dots ,((a_2,0,0),0,0), ((a_1,0,1),0,0),1] \in \mathcal{B}_{\lambda_q}\end{split}\end{equation*} (between neighbouring simples, here $a_1$ has to be strictly greater than $1$ in the second case) or 
\begin{equation*}\begin{split}&[((a_q,0,0),0,0), \dots,  ((a_{i+1},1,0),0,0),((p-a_i,0,0),1,0),\dots, ((a_1,0,0),0,0),0] ,\\&[((a_q,0,0),0,0), \dots,  ((a_{i+1},0,1),0,0),((p-a_i,0,0),1,0),\dots, ((a_1,0,0),0,0),0]  \in \mathcal{B}_{\lambda_q}\end{split}\end{equation*} 
for some $i$, between simples $L(\mu)$ and $L(\nu)$ with $m_\mu-m_\nu = 2(p-a_i)p^{i-1}$ in the first and $m_\mu-m_\nu =-2a_i p^{i-1}$ in the second case, which again needs $a_{i+1} >1$.

\section{Example.}\label{example}

The algebra $\lambda_q$ is isomorphic to the Yoneda extension
algebra of a block of Schur algebra $S(2,r)$ with $p^q$ simple
modules. Here, and throughout this paper, by the Yoneda extension
algebra of an abelian category, we mean the Yoneda extension algebra
of a complete set of nonisomorphic simple objects in that category;
by the Yoneda extension algebra of an algebra, we mean the Yoneda
extension algebra of the category of finite dimensional modules for
that algebra. To place our feet on the ground, let us describe an
example of such an algebra.

Let $p=3$. Let $\bfy_2$ denote the Yoneda extension algebra of a block
of a Schur algebra $S(2,r)$ with $9$ simple modules.  We label the $m$th simple module by $a_1,a_2$ if $(a_1-1)3+a_2 = m$. The algebra
$\bfy_2$ is isomorphic to $FQ/J$, where $Q$ is the quiver
$$
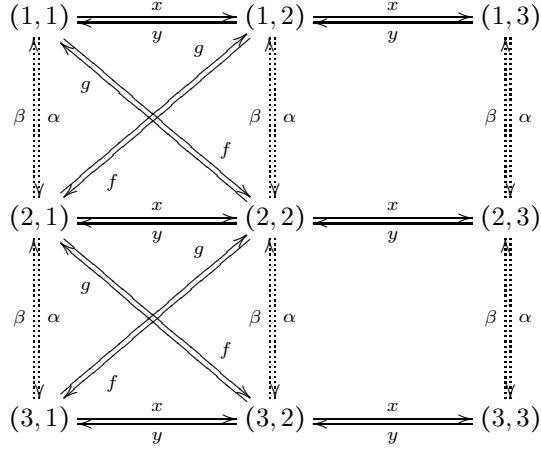
\xymatrix@=60pt{ (1,1) \ar[r]^x \ar@{..>}[d]^\alpha \ar[dr]^>>>>>>>{f} & \ar[dl]^>>>>>>>{f} \ar@<2pt>[l]^y (1,2) \ar@{..>}[d]^\alpha \ar[r]^x & \ar@<2pt>[l]^y (1,3) \ar@{..>}[d]^\alpha \\
(2,1) \ar@{..>}[d]^\alpha  \ar@{..>}@<2pt>[u]^\beta \ar[r]^x \ar[dr]^>>>>>>>{f} \ar@<2pt>[ur]^>>>>>>>{g} & \ar@<2pt>[ul]^>>>>>>>{g} \ar[dl]^>>>>>>>{f} \ar@<2pt>[l]^y (2,2) \ar@{..>}[d]^\alpha  \ar@{..>}@<2pt>[u]^\beta \ar[r]^x & \ar@<2pt>[l]^y (2,3) \ar@{..>}[d]^\alpha  \ar@{..>}@<2pt>[u]^\beta \\
(3,1) \ar@{..>}@<2pt>[u]^\beta \ar[r]^x \ar@<2pt>[ur]^>>>>>>>{g} & \ar@<2pt>[ul]^>>>>>>>{g} \ar@<2pt>[l]^y (3,2)  \ar@{..>}@<2pt>[u]^\beta \ar[r]^x & \ar@<2pt>[l]^y (3,3)  \ar@{..>}@<2pt>[u]^\beta
}$$
and $J$ is the ideal generated by the following relations
\begin{itemize}
\item $xy-yx, fg-gf, \alpha\beta-\beta \alpha$
\item $\begin{array}{cccc} x \beta - \beta x , & x \alpha - \alpha x &&\\
y \beta - \beta y ,&
y \alpha - \alpha y ,&&\\
f \beta - \beta f ,&
f \alpha - \alpha f ,&&\\
g \beta - \beta g ,&
g \alpha - \alpha g&& \\ yf-gx, & yg-fx, & gy-xf,& fy-xg ,\end{array}$
\item $yxe_{(i,1)}, gfe_{(1,i)}, \beta \alpha e_{(1,i)}, \beta f e_{(1,i)}, g \alpha e_{(1,i)}$ for $1 \leq i\leq 3$
\item $e_{(i,1)}gye_{(i+1,3)}$, $e_{(i+1,3)}xfe_{(i,1)}$, $e_{(i+1,1)}fye_{(i,3)}$, $e_{(i,3)}xge_{(i+1,1)}$ for $1 \leq i\leq 2$
\end{itemize}
where $e_{(i,j)}$ denotes the idempotent at vertex $(i,j)$.

The solid arrows have Yoneda degree $1$, whilst the dotted arrows
have Yoneda degree $3$. This presentation was obtained by direct computation, however,  the filtration of projectives by certain graded subquotients is directly visible from our construction.

For example, the projective indecomposable module indexed by vertex
$(1,1)$ has a basis given by paths

\begin{equation*}\begin{split} \mathfrak{B}:=
\{ &e, xe,x^2e,\\& fe,yfe, \\ &f^2e, xf^2e,  \\ &\alpha e,  x\alpha e,
x^2\alpha e,\\ & f\alpha e,yf\alpha e, \\ & \alpha^2 e,x\alpha^2
e,x^2\alpha^2 e\}\end{split}\end{equation*} where $e$ is the
idempotent corresponding to vertex $(1,1)$.

It's composition structure is given by 
$$
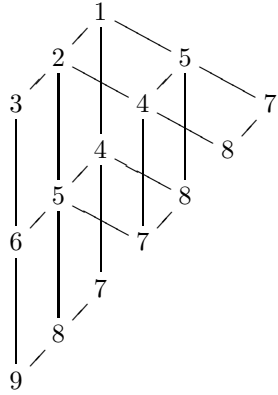
\xymatrix@=5pt{&& 1\ar@{-}[dl] \ar@{-}[drr] \ar@{-}[ddd] &&&&\\
& 2\ar@{-}[ddd]  \ar@{-}[dl]\ar@{-}[drr] &&& 5\ar@{-}[ddd] \ar@{-}[dl] \ar@{-}[drr] &&\\
3\ar@{-}[ddd]  &&& 4\ar@{-}[ddd] \ar@{-}[drr]  &&& 7\ar@{-}[dl] \\
&& 4 \ar@{-}[ddd] \ar@{-}[dl]\ar@{-}[drr] &&& 8 &\\
& 5\ar@{-}[ddd]  \ar@{-}[dl]\ar@{-}[drr] &&& 8\ar@{-}[dl] \\
6\ar@{-}[ddd]  &&& 7 &&&\\
&& 7 \ar@{-}[dl]&&&&\\
& 8 \ar@{-}[dl]&&&&&\\
9 &&&&&&
}$$
where as usual a line between two simples signifies that there is an extension between them and we have pictured components with Yoneda degree $d$ in row $d$. 
Writing the basis elements into the picture, we obtain
$$
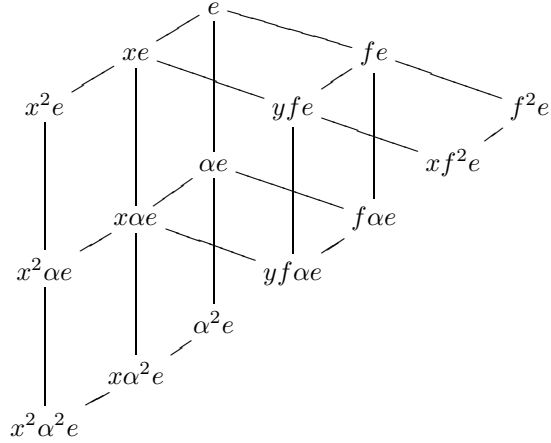
\xymatrix@=5pt{&& e\ar@{-}[dl] \ar@{-}[drr] \ar@{-}[ddd] &&&&\\
& xe\ar@{-}[ddd]  \ar@{-}[dl]\ar@{-}[drr] &&& fe\ar@{-}[ddd] \ar@{-}[dl] \ar@{-}[drr] &&\\
x^2e\ar@{-}[ddd]  &&& yfe\ar@{-}[ddd] \ar@{-}[drr]  &&& f^2e\ar@{-}[dl] \\
&& \alpha e \ar@{-}[ddd] \ar@{-}[dl]\ar@{-}[drr] &&& xf^2e &\\
& x\alpha e\ar@{-}[ddd]  \ar@{-}[dl]\ar@{-}[drr] &&& f\alpha e\ar@{-}[dl] \\
x^2\alpha e\ar@{-}[ddd]  &&& yf\alpha e &&&\\
&& \alpha^2 e \ar@{-}[dl]&&&&\\
& x\alpha^2 e \ar@{-}[dl]&&&&&\\
x^2\alpha^2 e &&&&&&
}$$
and the subquotients given by basis elements $\{e, xe,x^2e\}$,$\{\alpha e,  x\alpha e,
x^2\alpha e\}$ as well as $\{\alpha^2 e,x\alpha^2
e,x^2\alpha^2 e\}$ are each isomorphic to $\Omega e_1$, the subquotients given by basis elements $\{ fe,yfe\}$ and $\{f\alpha e,yf\alpha e\}$ are isomorphic to $\Theta^\sigma e_1$ and the subquotient given by basis elements $\{f^2e, xf^2e\}$ is isomorphic to $\Theta e_1$.

We now describe the monomial basis for this projective
indecomposable module. Let $\{ a,b,c \}$ be a basis of the
projective indecomposable $\Omega$-module indexed by vertex $1$, i.e.\ 
$$ a=(1,0,0), b=(1,1,0), c= (1,2,0) \in \mathcal{B}_\Omega,$$
lying in degrees $0,1,2$ respectively. 
Let $ \{\eta = (1,0,0), \theta= (1,1,0) \} \subset \mathcal{B}_\Theta$ and $\{ \xi = (2,0,0), \varsigma=(2,1,0) \} \subset \mathcal{B}_\Theta$ be the basis elements of $\Theta  e_1$ and
$\Theta^\sigma e_1$ respectively, which are not annihilated by the idempotent at the vertex $1$ from the right. The monomial
basis $\mathcal{B}_2$ is the set
\begin{equation*}\begin{split}
\{ &[(\gamma, 0,h),(\gamma', 0,h'),l],[(\gamma, 0,h),(\delta,n,h'),l],\\&
[(\delta, n,h),(\gamma,0, h'),l],[(\delta,n,h),(\delta', n',h'),l]  \}\end{split}\end{equation*} 
where $h,h',l \in \NN, n,n' \in
\NN_{>0}, \gamma \in \{a,b,c\}$ and $\delta \in \{ \eta, \theta\}$
or $\{ \xi, \varsigma\} $ depending on whether $n$ is even or odd.

The elements of height zero in here are
$$ \{  [(\gamma, 0,0),(\gamma',0,|\gamma|), p|\gamma| +|\gamma'|], [(\gamma, 0,0),(\delta,n, |\gamma|-n), |\delta| +p|\gamma|-pn  ]  \}$$
where $\gamma,\delta$ are before and $n=1, \dots, |\gamma|$.
So the full set is
\begin{equation*}\begin{split}
\mathfrak{B}':=\{&[(a,0,0),(a,0,0),0],[(a,0,0),(b,0,0),1], [(a,0,0),(c,0,0),2], \\
                                &  [(b,0,0),( \xi,1,0),0],[(b,0,0),( \varsigma,1,0),1],\\ &[(c,0,0),( \eta,2,0),0],[(c,0,0),( \theta,2,0),1]   \\
                                 &  [(b,0,0),(a,0,1),3],[(b,0,0),(b,0,1),4], [(b,0,0),(c,0,1),5]\\
                                  & [(c,0,0), ( \xi,1,1), 3], [(c,0,1), ( \varsigma,1,1), 4],\\
                                   & [(c,0,0),(a,0,2),6],[(c,0,0),(b,0,2),7], [(c,0,0),(c,0,2),8] \}\end{split}\end{equation*}

We can identify the bases $\mathfrak{B}$ and $\mathfrak{B}'$ for
$\bfy_2$, as ordered sets. The Yoneda degree of an element is given by summing the $n$s (i.e.\ the middle entries in the tuples) and $l$, giving $0,1,2,1,2,2,3,3,4,5,4,5,6,7,8$ respectively.

\section{Homological duality.} \label{kellerkoszul}

There are a number of approaches to homological duality in
representation theory. Here we describe two. One is a general
approach via differential graded algebras, due to Keller. The other
is a special approach for Koszul algebras, which is easier to work
with explicitly.

\subsection{Grading conventions.}\label{gradings}

In order to fix our sign conventions, we now give a brief introduction to dg algebras and modules.
A differential graded vector space is a $\mathbb{Z}$-graded vector space $V = \oplus_k V^k$ with a graded endomorphism $d$ of degree $1$.
We write $|v|_k$ for the degree of a homogeneous element of $V$.
We assume $d$ can act both on the left and the right of $V$, with the convention $d(v) = (-1)^{|v|_k}(v)d$.
A differential graded algebra is a $\mathbb{Z}$-graded algebra $A = \oplus_k A^k$ with a differential $d$ such that
$$d(ab) = d(a).b + (-1)^{|a|_k}a.d(b),$$
or equivalently
$$(ab)d = a.(b)d + (-1)^{|b|_k}(a)d.b.$$
If $A$ is a differential graded algebra then a differential graded left $A$-module is a graded left $A$-module $M$ with differential $d$ such that
$$d(a.m) = d(a).m + (-1)^{|a|_k}a.d(m);$$
a differential graded right $A$-module is a graded right $A$-module $M$ with differential $d$ such that
$$d(m.a) = d(m).a + (-1)^{|m|_k}m.d(a).$$
If $A$ and $B$ are dg algebras then a dg $A$-$B$-bimodule is a graded $A$-$B$-bimodule with a differential
which is both a left dg $A$-module and a right dg $B$-module.

If $_AM_B$ and $_BN_C$ are dg bimodules where $A$, $B$, and $C$ are dg algebras, then $M \otimes_B N$ is a dg $A$-$C$-bimodule
with differential $$d(m \otimes n) = d(m) \otimes n + (-1)^{|m|_k} m \otimes d(n).$$

Speaking about morphisms of dg algebras and dg (bi-)modules we mean homogeneous morphisms.
However, if $_AM_B$ and $_AN_C$ are dg bimodules where $A$, $B$, and $C$ are dg algebras, then $\Hom_A(M, N)$, the space of \emph{all} $A$-module homomorphisms from $M$ to $N$,  is a dg $B$-$C$-bimodule with differential $$d(\phi) = d\circ\phi - (-1)^{|\phi|_k} \phi\circ d.$$

If $_AM$ is a left dg $A$-module,
then $\End_A(M)$ is a differential graded algebra which acts on the right of $M$, giving $M$ the structure of an $A$-$\End_A(M)$-bimodule,
the differential on $\End_A(M)$ being given by $(\phi)d = \phi\circ d - (-1)^{|\phi|_k} d \circ\phi$.
If $M_B$ is a right dg $A$-module,
then $\End_A(M)$ is a differential graded algebra which acts on the left of $M$, giving $M$ the structure of an $\End_B(M)$-$B$-bimodule,
the differential on $\End_B(M)$ being given by $d(\phi) =  d\circ\phi - (-1)^{|\phi|_k} \phi \circ d$.

A differential bi- (tri-)graded vector space is a vector space $V$  with a $\mathbb{Z}^2$- respectively $\mathbb{Z}^3$-grading whose coordinates we denote by $(j,k)$ respectively $(i,j,k)$ and an endomorphism $d$ of degree $(0,0,1)$, ie.e $d$ is homogeneous with respect to the $i,j$-gradings and has degree $1$ in the $k$-grading, which we will also denote the homological grading.
We denote by $\langle \cdot \rangle$ a shift by $1$ in the $j$-grading, meaning $(V
\langle n \rangle)^{j} = V^{j-n}$. Since we will often identify dg modules and complexes, we will stick to the complex convention of $[\cdot]$ being a shift to the left, i.e.\ $V[n]_k = V_{k+n}$. Altogether $$(V\langle n \rangle [m])^{ijk} = V^{i,j-n,k+m}.$$
All definitions above can be extended to the differential bi- (tri-)graded setting, defining differential bi- (tri-)graded algebras, differential bi- (tri-)graded (left and right) $A$-modules as well as bi- (tri-)graded  $A$-$B$-bimodules as bi- (tri-)graded algebras resp.~modules resp.~bimodules which are differential graded algebras resp.~modules resp.~bimodules with respect to the $k$-grading, i.e. with respect to an endomorphism of degree $(0,0,1)$. Speaking about morphisms of differential bi- (tri-)graded  algebras and differential bi- (tri-)graded  (bi-)modules we mean homogeneous morphisms with respect to all gradings.
Similarly to the above, homomorphism spaces taken between $A$-modules (rather than differential (bi-) trigraded $A$-modules) will carry a differential bi- (tri-)grading.

For a dg algebra $A$, we denote by $D_{dg}(A)$ the dg derived category of $A$, whose objects are dg $A$-modules and where morphisms are given by the localisation of the set of dg module morphisms with respect to quasi-isomorphisms (see \cite[Section 3.1, 3.2]{Ke}).

We denote by $\mathbb{H}$ the
cohomology functor, which takes a differential $k$-graded complex
$C$ to the $k$-graded vector space $\mathbb{H} C = H^{\bullet} C$.

\subsection{The $2$-category $\mathcal{T}$.} Let $\mathcal{T}$ denote the $2$-category given by the following data:
\begin{itemize}
\item objects of $\mathcal{T}$ are pairs $(A,M)$ where $A = \bigoplus A^k$ is a dg algebra and $M
= \bigoplus M^k$ is a dg $A$-$A$-bimodule;
\item a $1$-morphism  $(X,\phi_X)$ between two objects $(B,N)$ and $(A,M)$ in $\mathcal{T}$ is given by a dg bimodule ${}_AX_B$ together with a quasi-isomorphism  of dg-bimodules $\phi_X: X \otimes_B N \rightarrow M \otimes_A X$;
\item a $2$-morphism between two such $1$-morphisms given by $({}_AX_B, \phi_X)$ and $({}_AY_B, \phi_Y)$ respectively is a morphism $f:X \to Y$ of dg bimodules such that the diagram 
$$\xymatrix{
 X \otimes_B N \ar^{\phi_X}[rr] \ar^{f \otimes id}[d]&& M \otimes_A X \ar^{id\otimes f}[d]\\
Y\otimes_B N \ar^{\phi_Y}[rr] && M \otimes_A Y
}$$ commutes.
\end{itemize}

\begin{defn}
We define a {\bf $j$-graded
object} of $\mathcal{T}$ to be an object $(a,m)$ of $\mathcal{T}$,
where $a= \bigoplus a^{jk}$ is a differential bigraded algebra, and $m
= \bigoplus m^{jk}$ a differential bigraded $a$-$a$-bimodule, and
$a^{jk} = m^{jk} = 0$ for $j<0$. 
\end{defn}

\begin{defn}
We define a {\bf Keller object} of $\mathcal{T}$ to be an object
$(A,M)$ of $\mathcal{T}$, where 
${}_AM_A$ is in $A\perfl \cap \perfr A$ and the natural morphisms of dg algebras $A\to \End_{A}({}_AM,{}_AM)$ and $A\to \End_{A}(M_A, M_A)$ are quasi-isomorphisms.
\end{defn}

We call such $M$ a two-sided tilting complex.

\begin{defn}
We define a {\bf Rickard object} of $\mathcal{T}$ to be a Keller object
$(A,M)$ of $\mathcal{T}$, where we have a quasi-isomorphism of dg algebras $A \to \bbH A$ and $\bbH A$ is a finite-dimensional algebra of finite-global dimension.
\end{defn}

This terminology is motivated by Keller's \cite{Ke} and Rickard's \cite{Rickard} Morita theory for dg derived categories and derived categories of algebras respectively, which identify objects with the properties asked of $M$ in each case as objects inducing autoequivalences of the respective categories.



\begin{defn}\label{dgeqdef} Let $(A,M)$ and $(B,N)$ be objects of $\mathcal{T}$.
 A {\bf dg equivalence} between objects $(A,M)$ and
$(B,N)$ of $\mathcal{T}$ is a pair $({}_AX_B,\phi_X)$ such that
\begin{enumerate}[(i)]
\item\label{dgeqdef1} ${}_AX$ belongs to $A \perf$, ${}_AX$ generates $D_{dg}(A)$, and the natural map
$$B \rightarrow \End_A(X)$$ is a quasi-isomorphism of dg algebras;
\item\label{dgeqdef2} $\phi_X\,:\,X \otimes_B N \rightarrow M \otimes_A X$ is a quasi-isomorphism of dg $A$-$B$-bimodules.
\end{enumerate}
\end{defn} 

Note that Definition \ref{dgeqdef}\eqref{dgeqdef2} just implies that $({}_AX_B,\phi_X)$ is a $1$-morphism in $\mathcal{T}$. Observe also that $X$ is a dg $A$-$B$-bimodule which induces an equivalence between $D_{dg}(A)$ and $D_{dg}(B)$ via
$$\xymatrix{ D_{dg}(A) \ar@/^/[r]^{\Hom_A(X,-)} & \ar@/^/[l]^{X \otimes_B -} D_{dg}(B)},$$
by Keller \cite[Lemma 3.10]{Ke}. We then have a diagram
$$\xymatrix{
D_{dg}(B)\ar^{X \otimes_B -}[r]\ar^{N \otimes_B -}[d] & D_{dg}(A)\ar^{M \otimes_A-}[d]\\
D_{dg}(B)\ar^{X \otimes_B -}[r] & D_{dg}(A) }$$
which commutes up to natural isomorphism of functors $X \otimes_B N \otimes_B - \rightarrow M \otimes_A X \otimes_B - $.

If there is a dg equivalence between $(A,M)$ and $(B,N)$,
we write $(A,M) \gtrdot (B,N)$.

\begin{defn}\label{qisodef} Let $(A,M),(B,N)\in\mathcal{T}$.
We define a {\bf quasi-isomorphism} from $(B,N)$ to $(A,M)$ to be a pair $(\varphi_B,\varphi_N)$
where $\varphi_B:B \rightarrow A$ is quasi-isomorphism of dg algebras, and $\varphi_N: {}_BN_B \rightarrow {}_AM_A$ is a quasi-isomorphism of vector spaces compatible with the bimodule structure and $\varphi_B$ in the sense that $\varphi_N(b_1nb_2)= \varphi_B(b_1)\varphi_N(n)\varphi_B(b_2)$. Similarly, we define an {\bf isomorphism} from $(B,N)$ to $(A,M)$ by replacing the word quasi-isomorphism by isomorphism throughout.
\end{defn} 

Note that this can be interpreted as a particular $1$-morphism $(A,\phi_A)$ from $(B,N)$  to $(A,M)$ in $\mathcal{T}$ where the right $B$-module structure on $A$ is induced by $\varphi_B$ and defining $\phi_A: A\otimes_B N \to M \otimes_A A : a \otimes n \mapsto a\varphi_N(n) \otimes 1$. The compatibility conditions in Definition \ref{qisodef} ensure that this is well-defined.

\subsection{Keller equivalence.}\label{kellereq} Here we recall a general statement concerning extension algebras of simple modules viewed as dg algebras. 

Let $A$ be a dg algebra, such that there is  a quasi-isomorphism of dg algebras $A \to \bbH A$ and $\bbH A$ is finite dimensional and of finite global dimension. Let $S_1, \dots, S_d$ be a complete set of non-isomorphic simple $\bbH A$-modules. Let $P_l = \bigoplus_k P^k_l$ be a projective $A$-resolution of $S_l$, viewed as a dg $A$-module. Note that since $A$ is assumed to finite global dimension, we can choose $P_l$ to lie in $A \perfl$. Denote by $\cE (A)$ the dg algebra 
$$\cE (A):=\bigoplus_{k,k'}\Hom_{A}(\bigoplus_{l=1}^dP^k_l,\bigoplus_{l=1}^dP^{k'}_l)$$
with the dg algebra structure as described in Section \ref{gradings}. Viewing these projective resolutions as complexes, the dg algebra differential translates into the total differential on the Hom double complex and therefore $\mathbb{H}\cE(A) = \Ext^\bullet(\bigoplus_{l=1}^d S_l, \bigoplus_{l=1}^d S_l).$
Then $P = \bigoplus_l P_l$ is a dg $A$-$\cE(A)$-bimodule.
There are mutually inverse equivalences
$$\xymatrix{ D_{dg}(A) \ar@/^/[r]^{\Hom_A(P,-)} & \ar@/^/[l]^{P \otimes_{\cE (A)}-} D_{dg}(\cE (A))},$$
by \cite[Lemma 3.10]{Ke}. Since $P$ is in $A \perfl$, we have a natural isomorphism of functors
$$\Hom_A(P,-) \cong \Hom_A(P,A) \otimes_A -.$$

Note that if $(A,M) \in \mathcal{T}$ is a Rickard object, the assumptions on $A$ in the above paragraph are automatically satisfied.

\subsection{The operator $\bbE$.} Here we extend the duality described in the previous section to objects in $\mathcal{T}$ by defining an operator $\bbE$ on $\mathcal{T}$. 

\begin{defn} Let $(A,M)$ be a Rickard object in $\mathcal{T}$ and $P$ be defined as in Section \ref{kellereq}.
Denote by $\cE(M) $ the dg $\cE(A)$-$\cE(A)$-bimodule
$$\cE(M) = \Hom_A(P,A) \otimes_A M \otimes_A P.$$  Furthermore, we define $$\bbE(A,M)
:=(\cE(A), \cE(M)).$$
\end{defn}

Note that $\cE(M) \cong \Hom_A(P,M) \otimes_A P$.

\begin{lem} \label{extsim} Let $(A,M)$ be a Rickard object in $\mathcal{T}$. The bimodule $P$, defined as in Section \ref{kellereq}, can be extended to a dg equivalence $(P,\phi_P)$ between $(A,M)$ and $\bbE(A,M)$.
\end{lem}
\proof We need to define a quasi-isomorphism of dg bimodules $\phi_P:P \otimes_{\cE(A)}
\cE(M) \rightarrow M \otimes_{A} P$, or in other words a
quasi-isomorphism
$$P \otimes_{\cE(A)} \Hom_A(P,A) \otimes_A M \otimes_A P \rightarrow M \otimes_A P.$$
We obtain this via the natural morphism of dg $A$-$A$-bimodules
$$P \otimes_{\cE(A)} \Hom_A(P,A) \rightarrow A,$$
which is an isomorphism since $P$ is a progenerator for
$A$.
\endproof

In the setup of Lemma \ref{extsim}, we thus have a natural isomorphism of functors making the following diagram
commute:
$$\xymatrix{D_{dg}(A) \ar^{\Hom_{A}(P,
-)}[rr] \ar^{M\otimes_A-}[d]&& D_{dg}(\cE
(A))\ar^{\cE(M)\otimes_{\cE(A)}-}[d]\\D_{dg}(A)\ar^{\Hom_{A}(P,
-)}[rr] && D_{dg}(\cE (A))}.$$

\begin{lem} \label{cancel}  Let $(A,M)$ be a Rickard object in $\mathcal{T}$ and $P$ defined as in Section \ref{kellereq}. We have an isomorphism of dg $\cE (A)$-$\cE (A)$-bimodules
$$\cE(M)^{\otimes_{\cE(A)} r} \rightarrow \Hom_A(P_A, M^{\otimes_A r} \otimes P_A).$$
\end{lem}
\proof If we write out $\cE(M)^{\otimes_{\cE(A)} r}$ in full, the internal occurences of the term $P \otimes_{\cE(A)} \Hom_A(P,A)$ cancel,
thanks to the isomorphism $P \otimes_{\cE(A)} \Hom_A(P,A)
\rightarrow A$ of dg $A$-$A$-bimodules, leaving us with
$$\Hom_A(P_A, A) \otimes_A M^{\otimes_A r} \otimes_A P_A \cong \Hom_A(P_A, M^{\otimes_A r} \otimes_A P_A).$$
\endproof

\subsection{Koszul duality.}\label{Koszulgen} In this subsection we summarise the results we need from \cite{BGS}. We assume $a = \bigoplus a^j$ is a graded
finite dimensional algebra with finite global dimension, with $a^0$
semisimple, and $a^j = 0$ for $j<0$. We further assume that $a$ is generated in degree $1$ and quadratic. Then
$$a = \mathbb{T}_{a^0}(a^1)/R,$$
where $R$ is the kernel of the multiplication map
$$a^1 \otimes_{a^0} a^1 \rightarrow a^2.$$
We write $x^* = \Hom({}_{a^0}x,{}_{a^0}a^0)$, and ${}^*x = \Hom(x_{a^0},
a^0_{a^0})$  if $x$ is a left/right $a^0$-module. If  $x$ is
concentrated in degree $j$, then by convention $x^*$ and ${}^*x$ are
concentrated in degree $-j$. The quadratic dual algebra $a^!$ is
given by
$$a^! = \mathbb{T}_{a^0}(a^{1*})/a^{2*},$$
where $a^{2*}$ embeds in $a^{1*} \otimes_{a^0} a^{1*}$ via the dual
of the multiplication map, composed with the inverse of the natural
isomorphism
$$a^{1*} \otimes_{a^0} a^{1*} \rightarrow (a^1 \otimes_{a^0} a^1)^*.$$
We then have $a^{!j} = 0$ for $j>0$.

The Koszul complex is defined to be the differential graded
$a$-$a^!$-bimodule
$$K = a \otimes_{a^0} {^*(a^!)},$$
whose $k$-grading is given by $$K^k = a \otimes_{a^0}
{}^*(a^{!-k}) \cong \Hom(a^{!-k}_{a^0}, a_{a^0}),$$ and whose
differential is given by the composition $$\Hom(a^{!-k-1}_{a^0},
a_{a^0}) \rightarrow \Hom(a^{!-k-1}_{a^0} \otimes_{a^0} a^1, a_{a^0}
\otimes_{a^0} a^1) \rightarrow \Hom(a^{!-k}_{a^0}, a_{a^0}),$$
obtained from the multiplication map $$a \otimes_{a^0} a^1
\rightarrow a,$$ and the natural composition $$a^{!-k} = a^{!-k}
\otimes_{a^0} a^0 \rightarrow a^{!-k} \otimes_{a^0} a^{1*}
\otimes_{a^0} a^1 \rightarrow a^{!-k-1} \otimes_{a^0} a_1.$$
Note that the differential has $j$-degree zero.

The algebra $a$ is said to be Koszul precisely when $K \rightarrow
a^0$ is a projective resolution of ${}_aa^0$, in which case the map
$a^! \rightarrow \cE (a)$ is a quasi-isomorphism. There are mutually
inverse equivalences
$$\xymatrix{ D^b(a \gr) \ar@/^/[r]^{\Hom_a(K,-)} & \ar@/^/[l]^{K \otimes_{a^!}-} D^b(a^! \gr)},$$
Under Koszul duality, the irreducible $a$-module $a^0e$ corresponds
to the projective $a^!$-modules $a^!e$, for a primitive idempotent
$e \in a^0$; a shift $\langle j \rangle$ in $D^b(a \gr)$ corresponds
to a shift $\langle j \rangle[-j]$ in $D^b(a^! \gr)$.

{\bf Warning.} Our $j$-grading convention does not coincide with
that in the literature. We assume that $a$ is generated in degrees
$0$ and $1$ whilst $a^!$ is generated in degrees $0$ and $-1$.
Beilinson, Ginzburg and Soergel assume that $a$ is generated in
degrees $0$ and $1$ but that $a^!$ is also generated in degrees $0$
and $1$. We have changed conventions in order to obtain gradings
which behave well with respect to our operators.

\subsection{The operator ${}^!$} In this subsection, we generalise the notion of Koszul duality from the previous subsection to objects of $\mathcal{T}$.

\begin{defn}\label{koszulop}
Suppose that $(a,m)$ is a $j$-graded Rickard object of $\mathcal{T}$, such that $a$ is a Koszul algebra.
We then call the object $(a,m)$ a {\bf Koszul object} of $\mathcal{T}$.

Let $\tilde{K}$ denote a projective resolution, viewed as a dg $a$-$a^!$-bimodule, of the complex $K$ of graded $a$-$a^!$-bimodules introduced in Section \ref{Koszulgen}. We define  $m^!$ to be the differential bigraded
$a^!$-$a^!$-bimodule $$m^!:=\Hom_a(\tilde{K},a) \otimes_a m \otimes_a \tilde{K}.$$
\end{defn}

\begin{lem}\label{basicdualities} Let $(a,m)$ be a Koszul object of $\mathcal{T}$ and $\tilde{K}$ as in Definition \ref{koszulop}.
\begin{enumerate}[(i)]
\item The dg bimodule $\tilde{K}$ can be extended a dg equivalence $(\tilde K, \phi_{\tilde K})$ between $(a,m)$ and $(a^!,m^!)$.
\item We have a quasi-isomorphism from $(a^!,m^!)$ to $\mathbb{E}(a,m)$.
\end{enumerate}
\end{lem}
\proof The proof of (i) is analogous to the proof Lemma
\ref{extsim}. The action of $a^!$ induces a quasi-isomorphism $a^! \rightarrow
\cE(a)$, from which part (ii) follows.
\endproof

\section{Operators $\mathbb{O}$ and $\fO$.} \label{ops}

\subsection{The operator $\mathbb{O}$.}  
In this subsection, we introduce the operator $\mathbb{O}$, which we showed to control the rational representation theory of $GL_2(F)$ in \cite{MT2}, cf.\ Section \ref{GL2}.

Let $A = \bigoplus A^k$ be a dg algebra, and $M = \bigoplus M^k$ a
dg $A$-$A$-bimodule.
The tensor algebra $\bbT_A(M)$ is differential bigraded with
$k$-degree by the total $k$-grading on each component $M^{\otimes j}$, with $A$ in
$j$-degree zero, and $M$ in $j$-degree $1$.

Let $(a,m)$ be a $j$-graded object of $\mathcal{T}$.
The tensor algebra $\bbT_a(m)$ is a differential trigraded algebra,
where $a$ has $i$-degree zero and $m$ has $i$-degree $1$, and the
$j$-$k$-bigrading is given by the total $j$- and the total $k$-grading on each component $m^{\otimes i}$.

\begin{defn} 
Given a $j$-graded object $(a,m)$ of
$\mathcal{T}$, we define an operator
$$\mathbb{O}_{a,m} \circlearrowright \mathcal{T}$$
given by
$$\mathbb{O}_{a,m}(A,M) = (\bigoplus a^{jk} \otimes_F M^{\otimes_A j}, \bigoplus m^{jk} \otimes_F M^{\otimes_A j}).$$
The algebra structure on $\bigoplus a^{jk} \otimes_F M^{\otimes_A j}$
is the restriction of the algebra structure on the super tensor product of
algebras $a \otimes \bbT_{A}(M)$. The $k$-grading and differential on
the complex $\bigoplus a^{jk} \otimes M^{\otimes_A j}$ are defined to
be the total $k$-grading and total differential on the tensor
product of complexes. The bimodule structure, grading and
differential on $\bigoplus m^{jk} \otimes M^{\otimes_A j}$ are defined
likewise.\end{defn}

Explicit formulae for the multiplication and differential on $a \otimes \bbT_{A}(M)$ are given as follows: For $\gamma, \delta \in a, \theta_1, \dots, \theta_r, \zeta_1, \dots, \zeta_s \in M$, 

\begin{equation*}\begin{split}(\gamma \otimes_F( \theta_1\otimes_A \dots\otimes_A \theta_r))&(\delta \otimes_F(\zeta_1\otimes_A \dots \otimes_A \zeta_s)) \\& = (-1)^{|\delta|_k(\sum_{t=1}^{r} |\theta_t|_k)} \gamma \delta \otimes_F  (\theta_1\otimes_A \dots\otimes_A \theta_r\otimes\zeta_1\otimes_A \dots \otimes_A \zeta_s)  \end{split}\end{equation*}

and

\begin{equation*}\begin{split}d(\gamma \otimes_F( \theta_1\otimes_A &\dots\otimes_A \theta_r)) =  d(\gamma) \otimes_F ( \theta_1\otimes_A \dots\otimes_A \theta_r)\\&+(-1)^{|\gamma|_k}\sum_{u=1}^r (-1)^{\sum_{t=1}^{u-1} |\theta_t|_k}\gamma \otimes_F( \theta_1\otimes_A \dots \otimes_A d(\theta_u)\otimes_A \dots\otimes_A \theta_r).
 \end{split}\end{equation*}

Similar formulae describe the bimodule structure on $m\otimes \bbT_{A}(M)$ by taking $\gamma$ or $\delta$ in $m$.

We again remark that this can be phrased in a $2$-categorical language, making $\mathbb{O}_{a,m}$ into a $2$-endofunctor of $\mathcal{T}$.

 We sometimes write
$$\mathbb{O}_{a,m}(A,M) = (a(A,M), m(A,M)).$$
If $a$ and $b$ are differential bigraded algebras, and ${}_a x_b$ is
differential bigraded $a$-$b$-bimodule, concentrated in nonnegative $j$-degrees, then we have a differential
graded $a(A,M)$-$b(A,M)$-bimodule $$x(A,M) := \bigoplus_{j,k} x^{jk} \otimes
M^{\otimes_A j}.$$

\begin{lem} \label{tensor}
Let $a, b,c$ be a differential bigraded algebras, ${}_bx_a$ and ${}_ay_c$
differential bigraded modules, all concentrated in nonnegative $j$-degrees. Let  $(A, M)$ be an object of $\mathcal{T}$.
Then
$$x(A,M) \otimes_{a(A,M)} y(A,M) \cong (x \otimes_a y)(A,M)$$ as differential bigraded $b(A,M)$-$c(A,M)$-bimodules.
\end{lem}
\proof We define a map $x(A,M) \otimes_{a(A,M)} y(A,M) \rightarrow
(x \otimes_a y)(A,M)$ sending homogeneous elements $(\alpha \otimes
m_1 \otimes \cdots \otimes m_{j_1}) \otimes (\beta \otimes n_1
\otimes \cdots \otimes n_{j_2})$ (where $\alpha \in x^{j_1}, \beta
\in y^{j_2}$ and $m_1, \dots m_{j_1},n_1,\dots, n_{j_2} \in M$ to
$(\alpha \otimes \beta) \otimes ( m_1 \otimes \cdots \otimes m_{j_1}
\otimes n_1 \otimes \cdots \otimes n_{j_2})$. It is straightforward to check that this is an isomorphism
of differential graded  $b(A,M)$-$c(A,M)$-bimodules.
\endproof

Suppose $(A,M)$ is a Keller object of $\mathcal{T}$. Then
$M^{-1} = \Hom_A(M,A)$ is a two-sided tilting complex such that $M^{-1}
\otimes_A - \cong \Hom_A(M,-)$ induces an inverse equivalence to $M
\otimes_A -$. The natural map $q: M \otimes_A M^{-1} \rightarrow A$
which takes $m \otimes \xi$ to $\xi(m)$ represents the counit of
the adjunction $(M \otimes_A -, M^{-1} \otimes_A -)$, and is a
quasi-isomorphism of dg $A$-$A$-bimodules.

We previously assumed that $(a,m)$ was a $j$-graded object of
$\mathcal{T}$ such that $a^{jk} = m^{jk} = 0$ for $j<0$. If we
assume rather that $a^{jk} = m^{jk} = 0$ for $j>0$, we define
$$\mathbb{O}_{a,m}(A,M) = (\bigoplus a^{jk} \otimes_F (M^{-1 })^{\otimes_A -j}, \bigoplus m^{jk} \otimes_F (M^{-1})^{ \otimes_A -j})$$
for a Keller object $(A,M)\in \mathcal{T}$.
Given a differential bigraded $a$-module $x$, with components in
positive and negative $j$-degrees, we define $x(A,M)$ to be the
$a(A,M)$-module given by
$$x(A,M) = (\bigoplus_{j<0} x^{j\bullet} \otimes (M^{-1})^{\otimes_A -j}) \oplus (x^{0\bullet} \otimes A) \oplus
(\bigoplus_{j>0} x^{j\bullet} \otimes M^{\otimes_A j}),$$
where the action is defined via the action of $a$ on $x$, along with
the quasi-isomorphism $q$. It is straightforward to check that this is well-defined.

\begin{lem} \label{homo}
Let $c$ be a differential bigraded algebra, $x$ and $y$ are differential bigraded
$c$-modules, all concentrated in nonnegative $j$-degrees, and let $(A,M)$ be a Rickard
object of $\mathcal{T}$. Then we have a quasi-isomorphism of differential bigraded vector spaces
$$\Hom_c(x,y)(A,M) \rightarrow \Hom_{c(A,M)}(x(A,M), y(A,M)).$$
\end{lem}
\proof We established this in a previous paper \cite[Proof of
Theorem 13]{MT2}. In that paper we only consider the case where
$(c,x=y)$ is a Rickard object. However, exactly the same proof works
in this more general case.\endproof

\begin{lem}\label{RickardtoKeller}
Suppose $(a,m)$ is a $j$-graded Rickard object  and $(A,M)$ is a Rickard
object of $\mathcal{T}$. Then $\bbO_{a,m}(A,M)$ is  Keller object of $\mathcal{T}$. 
\end{lem}
\proof We need to check that $m(A,M)$ is in $a(A,M)\perfl \cap \perfr a(A,M)$, which follows directly from the definition of $\bbO_{a,m}$ and that the natural maps $a(A,M) \to \End_{a(A,M)}({}_{a(A,M)}m(A,M))$ and $a(A,M) \to \End_{a(A,M)}(m(A,M)_{a(A,M)})$ are quasi-isomorphisms of dg algebras.
The latter follows since we have quasi-isomorphisms of dg algebras
\begin{equation*}
\begin{split}
a(A,M)  &\to \Hom_a({}_am,{}_am)(A,M) \\
&\to  \Hom_{a(A,M)}({}_{a(A,M)}m(A,M), {}_{a(A,M)}m(A,M))
\end{split}
\end{equation*}
and similarly for the other side. Here the first quasi-isomorphism comes from the fact that $(a,m)$ is a Rickard object, and the second follows from Lemma \ref{homo}.
 \endproof

\subsection{The operator $\fO$.} In this subsection, we introduce the operator $\fO$, whose relationship to $\mathbb{O}$ we will analyse in the next subsection, and whose main advantage is that it is well-behaved with respect to taking homology.

\begin{defn}
Let $\Gamma = \bigoplus \Gamma^{ijk}$ be a
differential trigraded algebra. We have an operator
$$\fO_\Gamma \circlearrowright \{ \Delta | \mbox{ $\Delta = \bigoplus \Delta^{jk}$ a differential bigraded algebra } \}$$
given by
\begin{equation}\label{fO}\fO_\Gamma(\Delta)^{ik} = \bigoplus_{j, k_1+k_2=k} \Gamma^{ijk_1} \otimes \Delta^{jk_2}.\end{equation}
The algebra structure and differential are obtained by restricting
the algebra structure and differential from $\Gamma \otimes \Delta$.
\end{defn}

If we forget the differential and the $k$-grading, the operator
$\fO_\Gamma$ is identical to the operator $\fO_\Gamma$ defined in the
introduction.

\begin{lem}\label{commute} We have $$\mathbb{H} \fO_{\Gamma} \cong \mathbb{H} \fO_{\mathbb{H} \Gamma} \cong \fO_{\mathbb{H} \Gamma}
\mathbb{H},$$ for a differential trigraded algebra $\Gamma$.
\end{lem} \proof Both isomorphisms are straightforward and based on the facts that
the tensor products in (\ref{fO}) are over $F$, so $\mathbb{H}  (\Gamma^{ijk_1} \otimes \Delta^{jk_2}) \cong \mathbb{H} ( \Gamma^{ijk_1}) \otimes \mathbb{H}( \Delta^{jk_2}),$ and that $\mathbb{H} \mathbb{H} =\mathbb{H} $. Using this we obtain

\begin{equation*}\begin{split}\mathbb{H}(\fO_\Gamma(\Delta))^{ik} &\cong \mathbb{H}(\bigoplus_{j, k_1+k_2=k} \Gamma^{ijk_1} \otimes \Delta^{jk_2})\\
&\cong \bigoplus_{j, k_1+k_2=k}\mathbb{H}(\Gamma^{ijk_1} ) \otimes  \mathbb{H}(\Delta^{jk_2}) .\end{split}\end{equation*}

Similarly, we have

\begin{equation*}\begin{split}\mathbb{H}( \fO_{\mathbb{H} \Gamma}(\Delta))^{ik} & \cong\mathbb{H}(\bigoplus_{j, k_1+k_2=k} \mathbb{H}(\Gamma^{ijk_1}) \otimes \Delta^{jk_2})\\
&\cong \bigoplus_{j, k_1+k_2=k}\mathbb{H}(\Gamma^{ijk_1} ) \otimes  \mathbb{H}(\Delta^{jk_2})
\end{split}\end{equation*}

and

\begin{equation*}\begin{split} (\fO_{\mathbb{H} \Gamma}\mathbb{H}(\Delta))^{ik} & \cong \bigoplus_{j, k_1+k_2=k}\mathbb{H}(\Gamma^{ijk_1} ) \otimes  \mathbb{H}(\Delta^{jk_2}),
\end{split}\end{equation*}

so they agree as graded vector spaces. Since all three multiplications are induced from $\mathbb{H} \Gamma\otimes\mathbb{H} \Delta$, the resulting differential bigraded algebras and therefore the operators are isomorphic.
\endproof

\subsection{Comparing $\mathbb{O}$ and $\fO$.} Here we describe relations
between the operators $\fO$ and $\mathbb{O}$. 
Throughout this subsection, let $(a,m)$ be a $j$-graded object of $\mathcal{T}$ and $(A,M) \in \mathcal{T}$.

Note that the algebra $\bbT_a(m)(A,M)$, formed with respect to the $j$-grading
on $\bbT_a(m)$, is a differential bigraded algebra, with
$$\bbT_a(m)(A,M)^{ik} = \oplus_{j} \bbT_a(m)^{ijk} \otimes M^{\otimes
j}.$$ The algebra $\bbT_{a(A,M)}(m(A,M))$ is a differential bigraded
algebra, with $$\bbT_{a(A,M)}(m(A,M))^{ik} =
(m(A,M)^{\otimes_{a(A,M)} i})^k.$$

We write $X^{i \diamond \bullet} \cong Y^{i \diamond \bullet}$ to signify that $X^{i jk} \cong Y^{ijk}$ for all $j,k$.

\begin{lem} \label{compareone}
\begin{enumerate}[(i)]
 \item\label{7i} We have an isomorphism of objects of $\mathcal{T}$
$$\bbO_{a,m}(A,M) =
(\fO_{\bbT_a(m)}(\bbT_A(M))^{0 \diamond \bullet},
\fO_{\bbT_a(m)}(\bbT_A(M))^{1 \diamond \bullet}),$$ where the $k$-grading
on the components of $\bbO_{a,m}(A,M)$ can be identified with the
$k$-grading on $\fO_{\bbT_a(m)}(\bbT_A(M))$.

\item \label{7ii} We have an isomorphism of differential bigraded algebras
$$\fO_{\bbT_a(m)}(\bbT_A(M)) \cong \bbT_a(m)(A,M).$$
\item\label{7iii} We have an isomorphism of differential bigraded algebras
$$\bbT_a(m)(A,M) \cong \bbT_{a(A,M)}(m(A,M)).$$
\end{enumerate}

\end{lem}
\proof The isomorphisms in \eqref{7i} and \eqref{7ii} follow directly from the definition, and the
isomorphism in \eqref{7iii} follows from Lemma \ref{tensor}, which gives $$m(A,M) \otimes_{a(A,M)} m(A,M) \cong (m \otimes_a m) (A,M).$$ \endproof

Suppose we are given $j$-graded objects $(a_i,m_i)$ in $\mathcal{T}$ for $1 \leq i \leq n$.
Let us define $(A_i,M_i)$ recursively via $(A_i,M_i) =
\bbO_{a_i,m_i}(A_{i-1},M_{i-1})$ and $(A_0,M_0) = (A,M)$.

\begin{lem} \label{comparetwo}
\begin{enumerate}[(i)]
 \item \label{8i} We have an algebra isomorphism
$$\bbT_{A_n}(M_n) \cong \fO_{\bbT_{a_n}(m_n)}
...\fO_{\bbT_{a_1}(m_1)}(\bbT_A(M)).$$
\item\label{8ii} We have an isomorphism
of objects of $\mathcal{T}$
$$\bbO_{a_n,m_n} ...\bbO_{a_1,m_1}(A,M) \cong$$
$$(\fO_{\bbT_{a_1}(m_1)}...\fO_{\bbT_{a_n}(m_n)}(\bbT_A(M))^{0\diamond \bullet},
\fO_{\bbT_{a_1}(m_1)}...\fO_{\bbT_{a_n}(m_n)}(\bbT_A(M))^{1\diamond \bullet}).$$
\end{enumerate}
\end{lem}

\proof The first part follows from Lemma \ref{compareone}\eqref{7ii},\eqref{7iii}  by induction. Part \eqref{8ii} then follows from Lemma \ref{compareone}\eqref{7i} and the first part.
\endproof

The application of Lemma \ref{comparetwo} which will be relevant for
our computation of $y$ is:

\begin{cor} \label{compare}
Let $(a,m)$ be a $j$-graded object of $\mathcal{T}$. Then we have an
isomorphism of differential bigraded algebras
$$\bbO_{F,0} \bbO^n_{a,m}(F,F) \cong \fO_F\fO_{\bbT_a(m)}^n(F[z]).$$
\end{cor}

\section{Operators $\bbO$ respect dualities.} \label{respect}

In this section, we analyse how the operators $\bbO$ introduced in Section \ref{ops} behave under the various incarnations of homological dualities that we have introduced in Section \ref{kellerkoszul}.

\subsection{Operators $\bbO$ and  dg equivalence.}
\begin{lem}\label{equivin}
Let $(A,M)$ and $(B,N)$ be objects of $\mathcal{T}$ such that $(A,M)
\gtrdot (B,N)$. Let $(a,m)$ be a $j$-graded Rickard object of
$\mathcal{T}$. Then $\mathbb{O}_{a,m}(A,M) \gtrdot
\mathbb{O}_{a,m}(B,N)$.
\end{lem}
\proof Let $X$ denote an $A$-$B$-bimodule inducing the dg
equivalence $(A,M) \gtrdot (B,N)$. We define $$x = a(A,M) \otimes_A
X$$ and claim that this is a dg $a(A,M)$-$a(B,N)$-bimodule which can be extended to a dg equivalence
$\mathbb{O}_{a,m}(A,M) \gtrdot \mathbb{O}_{a,m}(B,N)$. Clearly $x$
is a left differential bigraded $a(A,M)$-module. We have a sequence of canonical
homomorphisms of differential bigraded algebras
$$\begin{array}{ll}
a(B,N) & \cong B \otimes_B a(B,N) \\
& \cong \Hom_A(X, X) \otimes_B
a(B,N) \\
& \cong \Hom_A(X, X \otimes_B a(B,N)) \\
& \rightarrow \Hom_A(X, a(A,M) \otimes_A X) \\
& \cong \Hom_A(X, \Hom_{a(A,M)}(a(A,M),a(A,M) \otimes_A X)) \\
& \cong \Hom_{a(A,M)}(a(A,M) \otimes_A X, a(A,M) \otimes_A X) \\
& = \Hom_{a(A,M)}(x,x),
\end{array}$$
whose composition give $x$ the structure of a
$a(A,M)$-$a(B,N)$-bimodule. Here the third isomorphism is by virtue of $X$ being in $A \perfl$. Every term in this sequence is a
quasi-isomorphism, which means their composite is also a
quasi-isomorphism, as we require.

We have a sequence of quasi-isomorphisms of differential bigraded $a(A,M)$-$a(B,N)$-bimodules
$$\begin{array}{ll}
x \otimes_{a(B,N)} m(B,N) & = a(A,M) \otimes_{A} X \otimes_{a(B,N)} m(B,N) \\
& \rightarrow  (a \otimes_a m)(A,M) \otimes_A X \\
& \cong m(A,M) \otimes_A X  \\
& \cong m(A,M) \otimes_{a(A,M)} a(A,M) \otimes_A X  \\
& \cong m(A,M) \otimes_{a(A,M)} x
\end{array}$$
where the quasi-isomorphism is constructed analogously to the isomorphism of Lemma
\ref{tensor}. This defines our quasi-isomorphism $\phi_x$ needed to turn $(x,\phi_x)$ into a dg equvalence and completes the proof of Lemma \ref{equivin}.
\endproof

\begin{lem} \label{equivinagain} Let $(A,M)$ and $(B,N)$ be quasi-isomorphic Keller objects of $\mathcal{T}$.
Let $(a,m)$ be a $j$-graded object of $\mathcal{T}$. Then
$\mathbb{O}_{a,m}(A,M)$ and $\mathbb{O}_{a,m}(B,N)$ are
quasi-isomorphic objects of $\mathcal{T}$.
\end{lem}

\proof
Since by virtue of $(A,M)$ and $(B,N)$ being Keller objects, the endofunctors $M\otimes_A-, -\otimes_A M,N\otimes_B-, -\otimes_B N$ are exact, tensoring the quasi-isomorphism $f:{}_AM_A \to {}_BN_B$ together $i$ times produces the necessary quasi-isomorphisms  $f^i:{}_AM^{\otimes_A i}_A \to {}_BN^{\otimes_N i}_B$ in each $j$-graded component. It is immediate that this is compatible with the multiplicative resp.\ bimodule structures.\endproof

\subsection{Reversing gradings.}

Let us denote by $\mathbb{R}$ the sign reversing operator $\mathbb{R}$ on differential bigraded algebras,
which sends a differential bigraded algebra $a = \oplus_{j,k} a^{j,k}$ to itself with reverse grading, $\mathbb{R}a^{j,k} =
a^{-j,k}$. Likewise, if $m$ is a differential bigraded $a$-$a$-bimodule, we define $\mathbb{R} m$ by $\mathbb{R} m^{j,k} =
m^{-j,k}$.

\begin{lem}\label{minustoplus} Let $(a,m)$ and $(b,n)$ be graded
Rickard objects of $\mathcal{T}$, and $(A,M)$ a Rickard object of
$\mathcal{T}$. Then there is a quasi-isomorphism of objects of
$\mathcal{T}$, $$\bbO_{\mathbb{R}a, \mathbb{R}m}
\bbO_{b,n^{-1}}(A,M) \rightarrow \bbO_{a,m} \bbO_{b,n}(A,M).$$
\end{lem}
\proof For any $(B,N)\in \mathcal{T},$ we have an isomorphism
\begin{equation}\label{general} \mathbb{O}_{a,m}(B,N) \cong \mathbb{O}_{\mathbb{R}a, \mathbb{R}m}(B,N^{-1}),\end{equation}
since both sides are given by the object $$(\oplus_j a^j \otimes_F
N^{\otimes_A j}, \oplus_j m^j \otimes_F N^{\otimes_A j})$$ of
$\mathcal{T}$. 

 Note that $(b,n)$ being Rickard implies $(b,n^{-1})$ being Rickard, so, by Lemma \ref{RickardtoKeller},  both $\bbO_{b,n}(A,M)=(b(A,M), n(A,M))$ and $\bbO_{b,n^{-1}}(A,M) = (b(A,M), n^{-1}(A,M))$ are Keller objects. 

Since $(A,M)$ is Rickard, Lemma \ref{homo} gives us a quasi-isomorphism of dg $b(A,M)$-$b(A,M)$-bimodules
$$n^{-1}(A,M) \rightarrow \Hom_{b(A,M)}(n(A,M), b(A,M))$$
and  $(b(A,M), n(A,M))$ being a Keller object implies 
$$(b(A,M), n(A,M)^{-1}) = (b(A,M), \Hom_{b(A,M)}(n(A,M), b(A,M)))$$ is a Keller object.
We therefore have a quasi-isomorphism between Keller objects
\begin{equation}\label{Kellerqiso} (b(A,M), n^{-1}(A,M)) \rightarrow (b(A,M), n(A,M)^{-1})\end{equation}  in $\mathcal{T}$. The operator $\mathbb{O}_{\mathbb{R}a, \mathbb{R}m}$
respects quasi-isomorphisms between Keller objects by Lemma \ref{equivinagain}. Applying this to \eqref{Kellerqiso} and recalling \eqref{general} for the special case $(B,N) = \bbO_{b,n}(A,M)$, we obtain a
quasi-isomorphism
\begin{equation*}
\begin{split}
\bbO_{\mathbb{R}a, \mathbb{R}m} \bbO_{b,n^{-1}}(A,M)  &=\bbO_{\mathbb{R}a, \mathbb{R}m}(b(A,M), n^{-1}(A,M))\\
&\to \bbO_{\mathbb{R}a, \mathbb{R}m} (b(A,M), n(A,M)^{-1})\\
&\cong \bbO_{a,m}(b(A,M), n(A,M))\\
& = \bbO_{a,m} \bbO_{b,n}(A,M) 
\end{split}
\end{equation*}
in $\mathcal{T}$ as required.
\endproof

\subsection{Koszul duality for operators $\mathbb{O}$}
\label{Koszulops}

We now consider how Koszul duality behaves towards the operators
$\mathbb{O}$.

\begin{lem}\label{koszul}
Let $(A,M)$ be a Rickard object of $\mathcal{T}$. Let $(a,m)$ be a
Koszul object of $\mathcal{T}$, with Koszul dual $(a^!,m^!)$. Then
we have a dg equivalence $\mathbb{O}_{a,m}(A,M) \gtrdot 
\mathbb{O}_{a^!,m^!}(A,M)$.
\end{lem}
\proof Consider the differential bigraded $a(A,M)$-$a^!(A,M)$-bimodule $K(A,M)$. We show that $K(A,M)$ is can be extended to a dg equivalence $\mathbb{O}_{a,m}(A,M) \sim
\mathbb{O}_{a^!,m^!}(A,M)$ in $\mathcal{T}$. By Lemma \ref{tensor}, we have an
isomorphism of differential bigraded  $a(A,M)$-$a^!(A,M)$-bimodules
$$K(A,M) \cong a(A,M) \otimes_{a^0(A,M)} {^*(a^!)(A,M)}.$$
The action of $a^!$ on $K$ induces a quasi-isomorphism of differential bigraded algebras
$$a^! \rightarrow \Hom_a(K,K).$$
Consequently, using \cite[Lemma 15]{MT2}, we have a quasi-isomorphism of differential bigraded algebras
$$a^!(A,M) \rightarrow \Hom_a(K,K)(A,M)$$
which, when composed with the quasi-isomorphism $$\Hom_a(K,K)(A,M)
\rightarrow \Hom_{a(A,M)}(K(A,M),K(A,M))$$ of Lemma \ref{homo}, gives
us a quasi-isomorphism  of differential bigraded algebras
$$a^!(A,M) \rightarrow \Hom_{a(A,M)}(K(A,M),K(A,M))$$
which is induced by the action of $a^!(A,M)$ on $K(A,M)$.

The object $K(A,M)$ generates $D_{dg}(a(A,M)) \cong D_{dg}(a^!(A,M))$,
because $K(A,M)$ is quasi-isomorphic to $a^0(A,M) \cong a^0
\otimes_F A$, and $A$ generates $D_{dg}(A)$. All that remains for us
to do is to establish a quasi-isomorphism of differential bigraded
$a(A,M)$-$a^!(A,M)$-bimodules
$$K(A,M) \otimes_{a^!(A,M)} m^!(A,M) \rightarrow m(A,M) \otimes_{a(A,M)} K(A,M)$$
However, this follows from Lemma \ref{tensor}, and the fact that we
have a quasi-isomorphism of  differential bigraded $a$-$a^!$-bimodules
$$K \otimes_{a^!} m^! \rightarrow m \otimes_a K.$$
\endproof

\subsection{A quasi-isomorphism of operators} \label{qiops}

Here we show that, under good conditions, the operator $\bbE \bbO_{a,m}$ is quasi-isomorphic to
$\bbO_{a^!,m^!} \bbE$.

\begin{thm}\label{chain} Let $(A,M)$ be a Rickard object of $\mathcal{T}$.
Let $(a,m)$ be a be a Koszul object of $\mathcal{T}$, such that $ \bbO_{a,m}(A,M) $ is again a Rickard object. We have a
chain of dg equivalences
$$\bbE(\bbO_{a,m}(A,M))  \lessdot  \bbO_{a,m}(A,M) \gtrdot \bbO_{a^!,m^!}(A,M) \gtrdot \bbO_{a^!,m^!}(\bbE(A,M)). $$
\end{thm}

\proof This follows from Lemmas \ref{extsim}, \ref{equivin} and
\ref{koszul}.
\endproof

Theorem \ref{chain} implies we have a dg equivalence between the
objects $\bbE(\bbO_{a,m}(A,M))$ and $\bbO_{a^!,m^!}(\bbE(A,M))$ of
$\mathcal{T}$. We strengthen this as follows:

\begin{thm}\label{chainquasi}
The chain of equivalences in Theorem \ref{chain} lifts to a
quasi-isomorphism in $\mathcal{T}$ from $\bbO_{a^!,m^!}(\bbE(A,M))$ to
$\bbE(\bbO_{a,m}(A,M))$.
\end{thm}
\proof We need to prove that $\cE(a(A,M))$ is quasi-isomorphic to
$a^!(\bbE(A,M))$ as a dg algebra and that there is a compatible quasi-isomorphism between $\cE(m(A,M))$ and $m^!(\bbE(A,M))$ as in Definition \ref{qisodef}. We start with the first assertion. By definition,
$\cE(a(A,M)$ is obtained by the following recipe: first take a
projective resolution of the simple $a(A,M)$-modules, then take its
endomorphism ring. It is therefore enough to show that we can
compute $a^!(\bbE(A,M))$ by taking a projective resolution of the
simple $a(A,M)$-modules, and then taking the endomorphism ring.

We denote by $A^0$ a direct sum of a complete set of nonisomorphic
simple $A$-modules. Then $a^0 \otimes_F A^0$ is a direct sum of a
complete set of nonisomorphic simple $a(A,M)$-modules. Let
$P^\bullet_A$ denote a minimal projective resolution of ${}_AA^0$,
and $P^\bullet_a$ a minimal projective resolution of ${}_aa^0$. We
have $\cE(A)= \End(P^\bullet_A)$. A projective resolution of the
$a(A,M)$-module $a^0 \otimes A^0$ is given by $P_a(A,M) \otimes_A
P_A$. We can therefore write $\cE(a(A,M)) = \End(P_a(A,M) \otimes_A
P_A)$. We have a sequence of quasi-isomorphisms
$$\begin{array}{ll}
a^!(\cE(A), \cE(M))
& \cong \Hom_A(P_A, a^!(A,M) \otimes_A P_A) \\
& \mbox{\hspace{0.5cm} by Lemma \ref{cancel}}, \\
& \rightarrow \Hom_A(P_A, \cE(a)(A,M) \otimes_A P_A) \\
& \mbox{\hspace{0.5cm} by Koszul duality, and projectivity of $P_A$,} \\
& \cong \Hom_A(P_A, \Hom_a(P_a, P_a)(A,M) \otimes_A P_A) \\
& \mbox{\hspace{0.5cm} by definition of $\cE(a)$,} \\
& \cong \Hom_A(P_A, \Hom_{a(A,M)}(P_a(A,M), P_a(A,M)) \otimes_A P_A) \\
& \mbox{\hspace{0.5cm} by Lemma \ref{homo},} \\
& \cong \Hom_A(P_A, \Hom_{a(A,M)}(P_a(A,M), P_a(A,M) \otimes_A P_A) ) \\
& \mbox{\hspace{0.5cm} by projectivity and finite generation of $P_A$,} \\
& \cong \Hom_{a(A,M)}(P_a(A,M) \otimes_A
P_A, P_a(A,M) \otimes_A P_A) \\
& \mbox{\hspace{0.5cm} by adjunction,} \\
& \cong \End(P_a(A,M) \otimes_A P_A) \\
& = \cE(a(A,M)).
\end{array}$$

We likewise have a sequence of quasi-isomorphisms of dg vector spaces
$$\begin{array}{l}
m^!(\cE(A), \cE(M)) \\
\cong \Hom_A(P_A, m^!(A,M) \otimes_A P_A) \\
\rightarrow \Hom_A(P_A, \cE(m)(A,M) \otimes_A P_A) \\
\cong \Hom_A(P_A, \Hom_a(P_a, m \otimes_a P_a)(A,M) \otimes_A P_A) \\
\cong \Hom_A(P_A, \Hom_{a(A,M)}(P_a(A,M), m(A,M) \otimes_{a(A,M)} P_a(A,M)) \otimes_A P_A) \\
\cong \Hom_A(P_A, \Hom_{a(A,M)}(P_a(A,M), m(A,M) \otimes_{a(A,M)} P_a(A,M) \otimes_A P_A)) \\
\cong \Hom_{a(A,M)}(P_a(A,M) \otimes_A P_A, m(A,M) \otimes_{a(A,M)} P_a(A,M) \otimes_A P_A) \\
= \cE(m(A,M))
\end{array}$$
which are compatible with the bimodule structures in the sense of Definition \ref{qisodef}.
\endproof

In light of our interest in the object $\bbE(\bbO_{\bfc,\bft}^q(F,F))$
(see Section \ref{GL2}), we record the following corollary.

\begin{cor}\label{quasiiterated} Let $(a,m)$ be a be a Koszul object of $\mathcal{T}$, such that $\bbO_{a,m}^q (F,F)$ is a Rickard object for all $q \geq 0$. Then
 $\bbE ( \bbO_{a,m}^q (F,F))$ and
$\bbO_{a^!,m^!}^q (F,F)$ are quasi-isomorphic in $\mathcal{T}$.
\end{cor}

\proof This is proved by induction on $q$. The base step is given by
Lemma \ref{basicdualities} (ii), so we can inductively assume
$\bbE(\bbO_{a,m}^{q-1}(F,F))$ is quasi-isomorphic to
$\bbO_{a^!,m^!}^{q-1} (F,F)$. By Theorem \ref{chainquasi} and the
fact that $\bbO_{a,m}^{q-1}(F,F)$ is a Rickard object of
$\mathcal{T}$, the dg algebra $\bbE(\bbO_{a,m}^{q}(F,F))$ is
quasi-isomorphic to $\bbO_{a^!,m^!}\bbE(\bbO_{a,m}^{q-1}(F,F))$.
The latter is in turn, by Lemma \ref{equivinagain} and the
induction hypothesis, quasi-isomorphic to $\bbO_{a^!,m^!}^{q}
(F,F)$.\endproof

\section{Recollections of $GL_2$.} \label{GL2}

Let $Z$ denote the algebra given by the quiver
$$
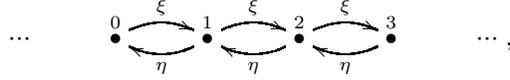
\xymatrix{ \cdots & \overset{0}{\bullet} \ar@/^/[r]^{\xi}
& \ar@/^/[l]^{\eta} \overset{1}{\bullet} \ar@/^/[r]^{\xi}
&\overset{2}{\bullet}\ar@/^/[l]^{\eta} \ar@/^/[r]^{\xi}
&\ar@/^/[l]^{\eta} \overset{3}{\bullet} &\cdots
 },$$
modulo relations $\xi^2 = \eta^2 = \xi \eta + \eta \xi = 0$. We will call this the \emph{zigzag algebra} and it or its various truncations show up in many guises in representation theory. It is therefore a very well-studied little infinite dimensional algebra whose projective indecomposable modules have a Loewy structure given by 
$$
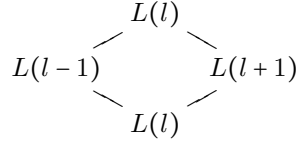
\xymatrix@=5pt{&L(l)\ar@{-}[dl]\ar@{-}[dr]&\\L(l-1)\ar@{-}[dr]&&L(l+1)\ar@{-}[dl]\\&L(l)&}
$$
where $L(l)$ denotes the simple module at vertex $l$, 
and is well-known to have
a number of interesting homological properties, all of which are easily checked by hand. For example, it is
quasi-hereditary, symmetric, and Koszul.

We denote by $\tau$ the algebra involution of $Z$ which sends vertex
$i$ to vertex $p-i$ and exchanges $\xi$ and $\eta$. Let $e_l$ denote
the idempotent of $Z$ corresponding to vertex $l \in \mathbb{Z}$.
Let
$$\bft = \sum_{1 \leq l \leq p, 0 \leq m \leq p-1} e_l Z e_m.$$ Then
$\bft$ admits a natural left action by the subquotient $\bfc$ of $Z$ given
by
$$\bfc =
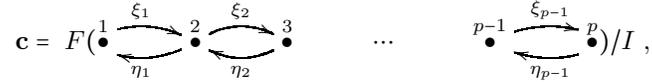
\xymatrix{F(\overset{1}{\bullet} \ar@/^/[r]^{\xi_1}
&\overset{2}{\bullet}\ar@/^/[l]^{\eta_1} \ar@/^/[r]^{\xi_2}
&\ar@/^/[l]^{\eta_2} \overset{3}{\bullet} &\cdots
&\overset{p-1}{\bullet}\ar@/^/[r]^{\xi_{p-1}} &
\ar@/^/[l]^{\eta_{p-1}} \overset{p}{\bullet})/I },$$ where
$I=(\xi_{l+1} \xi_{l}, \eta_{l} \eta_{l+1}, \xi_l \eta_l +
\eta_{l+1} \xi_{l+1} , \xi_{p-1} \eta_{p-1} \mid 1 \leq l \leq
p-2)$.  By symmetry, $\bft$ admits a right action by $\bfc$, if we twist
the regular right action by $\tau$. In this way, $\bft$ is naturally a
$\bfc$-$\bfc$-bimodule.
It is straightforward to see that the left restriction ${}_{\bfc}\bft$ of $\bft$ is a full characteristic tilting module for the quasi-hereditary algebra $\bfc$.
The natural homomorphism $\bfc \rightarrow \Hom({}_{\bfc} \bft,{}_{\bfc} \bft)$ defined by
the right action of $\bfc$ on $\bft$ is an isomorphism, implying that $\bfc$
is Ringel self-dual.

Note, in particular, that $\bfc$ and $\bft$ both graded by path length.
If we let $\tilde{\bft}$ denote a graded projective
resolution of $\bft$ as a $\bfc$-$\bfc$-bimodule, then $\tilde{\bft}$ is a
two-sided tilting complex, and $\tilde{\bft} \otimes_\bfc -$ induces a
self-equivalence of the derived category $D^b(\bfc)$ of $\bfc$.
Alternatively, we can (and will) interpret $\bfc$ as a differential bigraded algebra (with zero differential) with the $j$-grading given by path length, and by declaring $\bfc$ as concentrated in $k$-degree $0$. Then $\tilde{\bft}$ is a differential bigraded $\bfc$-$\bfc$-bimodule and $\bft$ is a differential bigraded $\bfc$-$\bfc$-bimodule concentrated in $k$-degree zero and with zero differential, making  $(\bfc,\tilde{\bft})$ and $(\bfc,\bft)$ into objects of $\mathcal{T}$.
In particular, since $\bfc$ as a quasi-hereditary algebra has finite global dimension, we have $\tilde{\bft} \in A \perfl$ and it being a $2$-sided tilting complex means that  $(\bfc,\tilde{\bft})$ is a Keller object. Moreover, since $\bfc$ equals its homology and is finite-dimensional and of finite global dimension, $(\bfc,\tilde{\bft})$ is  Rickard object of $\mathcal{T}$. The $j$-grading on both $\bfc$ and $\tilde{\bft}$ together with Koszulity of $\bfc$, which is easily checked by hand, now implies that  $(\bfc,\tilde{\bft})$ is even a Koszul object, which we will use extensively in the next section.

We will now recall some facts about the rational representation theory
of $G = GL_2(F)$. The category of polynomial representations of $G$
of degree $r$ is equivalent to the category $S(2,r) \ml$ of
representations of the Schur algebra $S(2,r)$ \cite{Green}. 
A block of $S(2,r)$ is Ringel self-dual if and only if it has $p^q$ simple modules \cite{EH}. In, \cite{MT1}.\cite{MT2} we developed a combinatorial way to describe these blocks, which we now describe.

The operator $\mathbb{O}_{\bfc,\bft}$ acts on $\mathcal{T}$, which containts the pair $(F,F)$ whose (dg) algebra is $F$ and whose (dg) bimodule is the regular bimodule $F$. Since both $\bfc$ and $\bft$ have zero differential, repeated application of $\mathbb{O}_{\bfc,\bft}$ to the pair $(F,F)$ produces an element of $\mathcal{T}$ where both components have zero differential and can hence regarded as an algebra and a bimodule respectively. The operator
$\mathbb{O}_{F,0}$ takes an object $(A,M)$ in $\mathcal{T}$ to the pair $(A,0)$, which we will identify with the dg algebra $A$.

We define $\bfb_q$ to be the category of modules over the
algebra $\mathbb{O}_{F,0} \mathbb{O}_{\bfc,\bft}^q(F,F)$. We have an
algebra homomorphism $\bfc \rightarrow F$ which sends a path in the
quiver to $1 \in F$ if it is the path of length zero based at $1$,
and $0 \in F$ otherwise. This algebra homomorphism lifts to a
morphism of operators $\mathbb{O}_{\bfc,\bft} \rightarrow
\mathbb{O}_{F,0}$. We have $\mathbb{O}_{F,0}^2 = \mathbb{O}_{F,0}$;
We thus have a natural sequence of operators
$$\mathbb{O}_{F,0} \leftarrow \mathbb{O}_{F,0} \mathbb{O}_{\bfc,\bft} \leftarrow \mathbb{O}_{F,0} \mathbb{O}_{\bfc,\bft}^2 \leftarrow ...,$$
which, if we apply each term to $(F,F)$ and take representations,
gives us a sequence of embeddings of abelian categories
$$\bfb_1 \rightarrow \bfb_2 \rightarrow \bfb_3 \rightarrow ...$$
We denote by $\bfb$ the union of these abelian categories. In a
previous paper \cite{MT2}, we have proved the following theorem:

\begin{thm}\cite[Corollary 21, Corollary 27]{MT2} \label{GL2theorem}
Every block of $G \ml$ is equivalent to $\bfb$. Every block of $S(2,r)
\ml$ whose number of isomorphism classes of simple objects is $p^q$
is equivalent to $\bfb_q$.
\end{thm}

To compute the Yoneda extension algebra $\bfY$ of $G \ml$, it is enough
to compute the Yoneda extension algebra $\bfy$ of the principal block
of $G \ml$ and, thanks to the above theorem, we can identify $\bfy$ with
the Yoneda extension algebra of $\bfb$. Since each embedding $\bfb_q
\rightarrow \bfb_{q+1}$ corresponds to taking an ideal in the poset of
irreducibles for $\bfb_{q+1}$, the theory of quasi-hereditary algebras \cite{CPS}
gives us a sequence of fully faithful embeddings of derived
categories 
$$D^b(\bfb_1) \rightarrow D^b(\bfb_2) \rightarrow D^b(\bfb_3) \rightarrow ...$$
Therefore, if we define $\bfy_q$ to be the Yoneda extension algebra of
$\bfb_q$, we have a sequence of algebra embeddings
$$\bfy_1 \rightarrow \bfy_2 \rightarrow \bfy_3 \rightarrow ...$$
whose union is equal to $\bfy$.

\section{Reduction.}\label{reduction}

\subsection{The pair $(\Omega,\Psi)$.}
Let $\bfc$ denote the quasi-hereditary algebra with $p$ irreducible
modules introduced in Section \ref{GL2}, and $\bft$ its tilting
bimodule and  $\tilde \bft$ a two-sided
tilting complex quasi-isomorphic to $\bft$. 

\begin{lem}\label{cqiso}
We have a quasi-isomorphism in $\mathcal{T}$ between 
$\bbO_{\bfc,\tilde \bft}^q(F,F)$ and  $\bbO_{\bfc,\bft}^q(F,F)$, in particular $\bbO_{\bfc,\tilde \bft}^q(F,F)$ is a Rickard object for any $q \geq 0$.
\end{lem}
\proof
In \cite[Lemma 22]{MT2}, we showed that for a quasi-hereditary Ringal self-dual algebra $A$, its tilting bimodule $T$ and a projective bimodule resolution $t$ of $T$, there is a quasi-isomorphism of differential bigraded algebras $\varphi_{1,1}:\bfc(A,t) \to \bfc(A,T)$. The proof relied only on the fact that $\bfc$ is concentrated in $j$-degrees $0,1,2$. As the same is true for $\bft$, we obtain a quasi-isomorphism $\varphi_{1,2}:\bft(A,t) \to \bft(A,T)$ which is compatible with $\varphi_{1,1}$ in the sense of Definition \ref{qisodef}.

Using this, we prove the lemma by induction on $q$, the case $q=1$ following by definition. Suppose the statement holds for $q-1$, i.e. we have a quasi-isomorphism $\bbO_{\bfc,\tilde \bft}^{q-1}(F,F) \to \bbO_{\bfc,\bft}^{q-1}(F,F)$, given by a quasi-isomorphisms $$\varphi_{q-1,1}:\bfc(\bbO_{\bfc,\tilde \bft}^{q-2}(F,F)) \to \bfc(\bbO_{\bfc,\bft}^{q-2}(F,F)) \hbox{ and } \varphi_{q-1,2}: \tilde\bft(\bbO_{\bfc,\tilde \bft}^{q-2}(F,F)) \to \bft(\bbO_{\bfc,\bft}^{q-2}(F,F)).$$ Since $\bft(\bbO_{\bfc,\bft}^{q-2}(F,F))$ is, by \cite[Proposition 20]{MT2}, the tilting bimodule for the quasi-hereditary algebra $\bfc(\bbO_{\bfc,\bft}^{q-2}(F,F))$, tensoring with it is exact on standard filtered modules (cf.\ \cite[Lemma 22]{MT2}), in particular on the tilting module itself. On the other hand, 
$\tilde\bft(\bbO_{\bfc,\tilde \bft}^{q-2}(F,F))$ is a two-sided tilting complex, hence tensoring with it is exact,  so we can tensor the latter quasi-isomorphism once to obtain a quasi-isomorphism
$$\tilde\bft(\bbO_{\bfc,\tilde \bft}^{q-2}(F,F))^{\otimes_{\bfc(\bbO_{\bfc,\tilde \bft}^{q-2}(F,F))} 2} \to \bft(\bbO_{\bfc,\bft}^{q-2}(F,F))^{\otimes_{\bfc(\bbO_{\bfc,\bft}^{q-2}(F,F))} 2}.$$ 
Since $\bfc$ is concentrated in degrees $0,1,2$, this suffices to obtain a quasi-isomorphism of differential bigraded algebras $\varphi_{q,1}: \bfc(\bbO_{\bfc,\tilde \bft}^{q-1}(F,F)) \to \bfc(\bbO_{\bfc,\bft}^{q-1}(F,F))$. In order to obtain  the necessary quasi-isomorphism $\varphi_{q,2}:\tilde\bft(\bbO_{\bfc,\tilde \bft}^{q-1}(F,F)) \to \bft(\bbO_{\bfc,\bft}^{q-1}(F,F))$ we compose  quasi-isomorphisms $$\tilde\bft(\bbO_{\bfc,\tilde \bft}^{q-1}(F,F)) \to \bft(\bbO_{\bfc,\tilde \bft}^{q-1}(F,F)) \to \bft(\bbO_{\bfc,\bft}^{q-1}(F,F))$$ where the second quasi-isomorphism is again constructed by virtue of $\bft$ being concentrated in degrees $0,1,2$ and the above considerations. Again $\varphi_{q,2}$ is compatible with the bimodule structures as required in Definition \ref{qisodef}, hence this shows that $\bbO_{\bfc,\tilde \bft}^q(F,F)$ and $\bbO_{\bfc,\bft}^q(F,F)$ are quasi-isomorphic in $\mathcal{T}$. 

To see that $\bbO_{\bfc,\tilde \bft}^q(F,F) = ( \bfc(\bbO_{\bfc,\tilde \bft}^{q-1}(F,F)),\tilde\bft(\bbO_{\bfc,\tilde \bft}^{q-1}(F,F))) $ is Rickard, we note that  it is Keller since $\tilde\bft(\bbO_{\bfc,\tilde \bft}^{q-1}(F,F))$ is a two-sided tilting complex. Furthermore $\varphi_{q,1}$ gives us a quasi-isomorphism between $\bfc(\bbO_{\bfc,\tilde \bft}^{q-1}(F,F))$ and its homology $\bfc(\bbO_{\bfc,\bft}^{q-1}(F,F))$, which is a finite-dimensional algebra of finite global dimension, as it is Morita equivalent to a block of a Schur algebra (cf.~Section \ref{GL2}). Therefore $\bbO_{\bfc,\tilde \bft}^q(F,F)$ is Rickard.
\endproof

We can hence apply our homological theory to the Koszul object $(\bfc,\tilde \bft) \in \mathcal{T}$ (cf.~Section \ref{GL2}) and thus define a
differential bigraded $\bfc^!$-$\bfc^!$-bimodule $\bft^!$ according to Definition \ref{koszulop}. In light of Corollary \ref{quasiiterated},
we wish to compute (the homology of)  iterated applications
of the operator $\bbO_{\bfc^!,\bft^!}$ to the pair $(F,F)$. However, $\bfc^!$
is now negatively graded so to evaluate this we need to work with
the adjoint of $\bft^!$, which is given by the differential bigraded
$\bfc^!$-$\bfc^!$-bimodule
$$\bft^{!-1} = \Hom_{\bfc^!}(\bft^{!},
\bfc^!).$$ Since $\bft^!$ is in $\bfc^! \perfl$ (it is quasi-isomorphic to the tilting module for $\bfc^!$, and the latter, as a quasi-hereditary algebra, has finite global dimension), we have an isomorphism of functors
$$\bft^{!-1} \otimes_{\bfc^!} - \cong \Hom_{\bfc^!}(\bft^!,-).$$
For notational simplicity, we wish to describe algebras with positive
gradings rather than negative gradings. Therefore rather than
working with the negatively graded algebra $\bfc^!$, we work with the
positively graded algebra $\Omega$, which is canonically isomorphic
to $\bfc^!$ as an algebra, but whose $j^{th}$ homogeneous component is
equal to the $-j^{th}$ homogeneous component of $\bfc^!$. We define the
differential bigraded $\Omega$-$\Omega$-bimodule $\Psi$ to be the
differential bigraded $\bfc^!$-$\bfc^!$-bimodule $\bft^{!-1}$, whose
$(j,k)^{th}$ homogeneous component is equal to the $(-j,k)^{th}$
homogeneous component of $\bfc^!$. Thus 

$$\Omega = \mathbb{R} \bfc^! \qquad\hbox{ and } \qquad\Psi = \mathbb{R} \bft^{!-1}.$$

Note that $$\Omega = \xymatrix{F(\overset{1}{\bullet}
\ar@/^/[r]^{x_1} &\overset{2}{\bullet}\ar@/^/[l]^{y_1}
\ar@/^/[r]^{x_2} &\ar@/^/[l]^{y_2} \overset{3}{\bullet} &\cdots
&\overset{p-1}{\bullet}\ar@/^/[r]^{x_{p-1}} & \ar@/^/[l]^{y_{p-1}}
\overset{p}{\bullet})/I^\perp },$$ where $I^\perp=(x_l y_l - y_{l+1}
x_{l+1} , y_1x_1 \mid 1 \leq l \leq p-2)$.

Since we consider $\Omega$ as the extension algebra of $\bfc$ (with the $j$-grading multiplied by $-1$), it is naturally a dg algebra concentrated in positive $k$-degrees with zero differential. The generators $x$ and $y$ are in $k$-degree 1.
For our computation of $\bfy$ we will use the following corollary of
Lemma \ref{minustoplus}:

\begin{cor} \label{minustopluscor} We have a quasi-isomorphim of differential graded algebras
$$\bbO_{F,0} \bbO_{\Omega,\Psi}^q(F,F) \to \bbO_{F,0} \bbO_{\bfc^!,\bft^!}^q (F,F).$$
\end{cor}
\proof We prove a slightly stronger statement, namely that  $$\bbO_{\Omega,
\Psi}^q(F,F) \cong (\bfc(\bbO_{\bfc^!,\bft^!}^{q-1} (F,F)), \bft^{!-1}(\bbO_{\bfc^!,\bft^!}^{q-1} (F,F))) = \bbO_{ \bfc^!,  \bft^{!-1}}(\bbO_{\bfc^!,\bft^!}^{q-1} (F,F)) $$ from which the Corollary follows by applying $\bbO_{F,0}$.

We proceed by induction on $q$. The case $q=1$ follows from the definitions of $\mathbb{R} \bfc^! = \Omega$ and $\mathbb{R} \bft^{!-1} =\Psi$.
Now assume we have a quasi-isomorphism $$\bbO_{\Omega,
\Psi}^{q-1}(F,F) \to (\bfc(\bbO_{\bfc^!,\bft^!}^{q-2} (F,F)), \bft^{!-1}(\bbO_{\bfc^!,\bft^!}^{q-2} (F,F))).$$
Then
\begin{equation*}\begin{split}
\bbO_{\Omega,\Psi}^q(F,F)&= \bbO_{\mathbb{R} \bfc^!, \mathbb{R} \bft^{!-1}}(\bbO_{\Omega,\Psi}^{q-1}(F,F))
\\
&\to \bbO_{\mathbb{R} \bfc^!, \mathbb{R} \bft^{!-1}}(\bfc(\bbO_{\bfc^!,\bft^!}^{q-2} (F,F)), \bft^{!-1}(\bbO_{\bfc^!,\bft^!}^{q-2} (F,F)))\\
&\to \bbO_{ \bfc^!,  \bft^{!-1}}(\bfc(\bbO_{\bfc^!,\bft^!}^{q-2} (F,F)), \bft^!(\bbO_{\bfc^!,\bft^!}^{q-2} (F,F)))\\
& = \bbO_{ \bfc^!,  \bft^{!-1}}(\bbO_{\bfc^!,\bft^!}^{q-1} (F,F))
\end{split}\end{equation*}
where the first equality is just by definition, the second isomorphism is by the inductive assumption and Lemma \ref{equivinagain}, the third is the quasi-isomorphism given in Lemma \ref{minustoplus} and the final equality is just the definition.
\endproof

\subsection{A chain of isomorphisms.}
Let us recapitulate what we have done so far towards our goal of computing the extension algebra $\bfy_q$ of a block of $GL_2(F)$ with $p^q$ simples modules. We have the following sequence of isomorphisms

$$\begin{array}{lll}
\bfy_q  & \cong \mathbb{H} \bbO_{F,0} \bbE \bbO_{\bfc,\bft}^q(F,F) & \mbox{Theorem \ref{GL2theorem}} \\
 & \cong \mathbb{H} \bbO_{F,0} \bbE \bbO_{\bfc,\tilde\bft}^q(F,F) & \mbox{Lemma \ref{cqiso}} \\
& \cong \mathbb{H} \bbO_{F,0} \bbO_{\bfc^!,\bft^!}^q \bbE(F,F) & \mbox{Corollary \ref{quasiiterated}} \\
& \cong \mathbb{H} \bbO_{F,0} \bbO_{\bfc^!,\bft^!}^q(F,F) & \mbox{$\bbE(F,F)=(F,F)$} \\
& \cong \mathbb{H} \bbO_{F,0} \bbO_{\Omega, \Psi}^q(F,F) & \mbox{Corollary \ref{minustopluscor}}\\
& \cong \mathbb{H} \fO_F \fO_{\bbT_{\Omega}(\Psi)}^q(F[z]) & \mbox{Corollary \ref{compare}} \\
& \cong \fO_{\mathbb{H} F} \fO_{\mathbb{H}\bbT_{\Omega}(\Psi)}^q \mathbb{H}(F[z]) & \mbox{Lemma \ref{commute}} \\
& \cong \fO_{F} \fO_{\mathbb{H}\bbT_{\Omega}(\Psi)}^q (F[z]) & \mathbb{H} F = F,\mathbb{H}(F[z]) =F[z].  \mbox{}
\end{array}$$

Recalling our definition of $\lambda_q =  \fO_F \fO_\Lambda^q(F[z])$ and our goal of proving an isomorphism of $\bfy_q$ and $\lambda_q$, we therefore focus our attention on
%
%
on proving the following result:

\begin{thm} \label{Lambdacalculation} We have an isomorphism of
$ijk$-trigraded algebras
$$\mathbb{H}(\bbT_{\Omega}(\Psi)) \cong \Lambda.$$
\end{thm}

This will take up most of the remainder of this article, and is largely independent of the preceding subsections, whose aim it was to establish the above sequence of isomorphisms.

\section{An explicit dg algebra.} \label{explicit}


\subsection{The homology $\mathbb{H}(\Psi)$.}
This subsection is concerned with verifying the statement of Theorem \ref{Lambdacalculation} for $i$-degree one.

In the following, we will write differential $k$-graded modules as complexes by specifying the $0$-th homological position and interpreting $\cdots \to M_s \to M_{s+1} \to \dots$ as $\bigoplus M_s[-s]$.

{\bf Warning:} The modules $M_s$ may carry their own $k$-grading, and will generally do so whenever we look at modules over $\Psi$, which, unlike $\bfc$, is not an algebra concentrated in $k$-degree zero. 

In order to simplify our computations and to use results from the theory of Koszul duality in their original setting, we will therefore initially work in the (bounded) derived category of $j$-graded $\Omega$-modules (e.g.~in Lemma \ref{csimples}, Lemma \ref{specitilt}, Lemma \ref{ses},  Lemma \ref{Lambda1}\eqref{26i},\eqref{26ii}), ignoring any internal $k$-grading of the algebra and then establishing the $k$-grading on the $j$-graded bimodules subsequently.

In order to compute the homology of $\Psi$, which we need for our computation of $\mathbb{H}(\bbT_{\Omega}(\Psi))$, we first compute the functor $\bft^{-1} \otimes -$ on various $\bfc$-modules.

The algebra $\bfc$ has irreducible modules $L(l)$ indexed by integers
$1 \leq l \leq p$. The algebra has standard modules $\Delta(l)$
which have top $L(l)$ and socle $L(l-1)$ in case $2 \leq l \leq p$,
and $\Delta(1) = L(1)$. The algebra has costandard modules
$\nabla(l)$ which have socle $L(l)$ and top $L(l-1)$ in case $2 \leq
l \leq p$, and $\nabla(1) = L(1)$. We define $L(l)$, $\Delta(l)$ and
$\nabla(l)$ as $j$-graded modules by insisting their tops are
concentrated in degree $0$.

\begin{lem}\label{csimples}
We have exact triangles in the derived category of $j$-graded $\bfc$-modules as follows:
$$L(p-l) \rightarrow \bft^{-1} \otimes_\bfc L(l) \rightarrow L(p)\langle -l \rangle[1-l] \rightsquigarrow,$$
for $1 \leq i < p$, and
$$0 \rightarrow \bft^{-1} \otimes_\bfc L(p) \rightarrow L(p) \langle -p \rangle [1-p] \rightsquigarrow.$$
\end{lem}
\proof We can identify $\bft^{-1} \otimes_{\bfc} -$ with $\Hom_{\bfc}(\bft,-)$, which by
Ringel duality takes costandard modules $\nabla$ to standard modules
$\Delta$. Easy computations establish the following formulae:
$$\bft^{-1} \otimes_\bfc \nabla(l) = \Delta(p-l+1),$$ for $2 \leq
l \leq p$ and
$$\bft^{-1} \otimes_\bfc \nabla(1) = \Delta(p) \langle -1 \rangle.$$

Since we have a quasi-isomorphism
$$L(l) \cong (\nabla (l) \langle -1\rangle \rightarrow \nabla (l-1) \langle -2\rangle \rightarrow \dots \rightarrow \nabla (2)\langle -(l-1)\rangle \rightarrow \nabla (1)\langle -(l-1)\rangle$$
we obtain
$$\bft^{-1} \otimes_\bfc L(l) = (\Delta(p-l+1) \langle -1 \rangle \rightarrow \Delta(p-l+2) \langle -2 \rangle \rightarrow ... \rightarrow  \Delta(p) \langle -l \rangle)$$
for $1 \leq l \leq p$, where in both complexes the leftmost terms lie in homological position
$0$. For $1 \leq l \leq p-1$, this is exact except in homological position
$0$, where the kernel is $L(p-l)$ and $k$-degree $l-1$, where the
cokernel is $L(p) \langle -l \rangle$. For $l=p$, this is exact
except in homological position $p-1$, where the cokernel is $L(p)\langle -p
\rangle$. Noting that a simple $L$ in homological position $t$ is the same as $L[-t]$, we obtain the statement of the lemma.
\endproof

We now compute a one-sided tilting complex which is quasi-isomorphic
to $\Psi$. The advantage of this complex over $\Psi$ is that its
structure is extremely explicit, making computations possible.

Consider the differential $jk$-bigraded $\Omega$-module which,
written as a two-term complex $W = (W^1 \overset{d_W}{\rightarrow} W^0)$, is given by
$$\Omega e_p \otimes_{F} e_pFQe  \rightarrow \Omega e$$ with $\Omega e$ in homological position zero,
where $Q$ is the type $A_p$ subquiver of the quiver of $\Omega$
generated by arrows $x$, where $FQ$ is the corresponding hereditary
subalgebra of $\Omega$, where $e = \sum_{1 \leq l \leq p-1}e_l$ is
the sum over the idempotents $e_l$ at $l$ and where the differential
on the complex is given by the algebra product. 
In other words, $W$ is given by
$$\bigoplus_{l=1}^{p-1}(\Omega e_{p} \langle l \rangle \to \Omega e_{p-l}),$$
where the maps are given by right multiplication with $x^l$.
Consider also the differential $jk$-bigraded $\Omega$-module which,
written as a two term complex, is given by
$$\Omega e_p \langle p\rangle \rightarrow 0$$ with $0$ in homological position $0$.

Let $X$ denote the sum of these two-term complexes.

\begin{lem} \label{specitilt}
We have a quasi-isomorphism between complexes of $j$-graded $\Omega$-modules $X$ and $\Psi$.
\end{lem}
\proof Under Koszul duality
an irreducible $\bfc$-module corresponds
to a projective $\bfc^!$-module, hence a projective $\Omega$-module.
A shift $\langle j \rangle$ for $\bfc$-modules corresponds to a shift $\langle j
\rangle[-j]$ for $\bfc^!$-modules, therefore a shift $\langle -j \rangle[-j]$ for $\Omega$-modules
(again the homological grading in the derived category of $\Omega$-modules does not coincide
with the $k$-grading since in the first $_\Omega\Omega$ is concentrated in a
single degree whilst in the second $_\Omega\Omega$ is concentrated in many degrees).
Lemma \ref{csimples} therefore gives us exact triangles
in the derived category of $j$-graded $\Omega$-modules as follows:
$$\Omega e_{p-l} \rightarrow \Psi \otimes_{\Omega} \Omega e_{l} \rightarrow \Omega e_{p} \langle l \rangle[1] \rightsquigarrow,$$
for $1 \leq l < p$, and
$$0 \rightarrow \Psi \otimes_{\Omega} \Omega e_{p} \rightarrow \Omega e_{p} \langle p \rangle[1] \rightsquigarrow.$$
Here $\Omega e_{l}$ denotes the projective indecomposable
$\Omega$-module indexed by $l$. We therefore have exact triangles in the derived category of  $j$-graded $\Omega$-modules
$$\Omega e_{p} \langle l \rangle \rightarrow \Omega e_{p-l} \rightarrow \Psi \otimes_{\Omega} \Omega e_{l}  \rightsquigarrow,$$
for $1 \leq l < p$. Morphisms in this derived category between
projective objects lift to morphisms of chain complexes. The object
$\Psi \otimes_{\Omega} \Omega e_{l}$ is indecomposable. Up to a scalar,
there is a unique graded homomorphism $\phi$ from $\Omega e_{p}
\langle l \rangle$ to $\Omega e_{p-l}$. This implies we have an
isomorphism in the derived category of $j$-graded $\Omega$-modules
$$\Psi \otimes_{\Omega} \Omega e_{l} \cong (\Omega e_{p} \langle l \rangle \overset{\phi}{\rightarrow} \Omega e_{p-l}),$$
or to put it another way, an isomorphism
$$\Psi \otimes_{\Omega} \Omega e_{l} \cong (\Omega e_p \otimes_{F} e_pFQe_{p-l} \rightarrow \Omega e_{p-l}),$$
where the term $\Omega e_{p-l}$ is concentrated in homological degree $0$,
and the differential is given by multiplication. We thus find $X$ is
quasi-isomorphic to $\Psi$ as a complex of $j$-graded $\Omega$-modules  as required.
\endproof

\begin{lem} \label{ses} We have an exact sequence of $j$-graded $\Omega$-modules,
$$0 \rightarrow \Omega e_l \langle p \rangle \rightarrow \Omega e_p \langle l \rangle \rightarrow \Omega e_{p-l} \rightarrow \Theta e_{p-l} \rightarrow 0,$$
for $1 \leq l \leq p-1$.
\end{lem}

\proof 
This is a direct computation visible from the Loewy structure of these small modules, as described in Section \ref{intro}.
%
%
%
\endproof

Taking the first two terms in the sum from $1$ to $p-1$ of the exact sequences in Lemma \ref{ses}, we obtain 
 $W^0 \langle p \rangle$ and $W^1$, respectively. We define the sum
of the maps between these terms to be
$$\delta: W^0 \langle p \rangle \rightarrow W^1.$$

Let $\gamma$ denote the automorphism of $\Omega$ that sends a homogeneous element $\omega$ to $(-1)^{|\omega|_k} \omega$. The following is an easy direct computation.

\begin{lem}\label{inner}
$\gamma \circlearrowleft \Omega$ is an inner automorphism that sends $\omega$ to $g \omega g^{-1}$, where $g = \sum (-1)^l e_l$.
\end{lem}

The next lemma establishes the aim of this subsection, by showing that the algebras $\Lambda$ and $\mathbb{H}(\Psi)$ are isomorphic in $i$-degree one.

\begin{lem}\label{Lambda1}
\begin{enumerate}[(i)]
\item\label{26i} $X$ is a $j$-graded tilting complex for $\Omega$.
\item\label{26ii} We have an isomorphism of $j$-graded algebras
$$\mathbb{H} \End_{\Omega}(X) \cong \Omega.$$
\item\label{26iii} As a $jk$-graded left $\Omega$-module, the kernel of the differential on $X$ is isomorphic to $\Omega
\langle p \rangle[1-p]$, and the cokernel of the differential on $X$ is
isomorphic to $\Theta^{\sigma \gamma}$.
\item\label{26iv} We have isomorphisms of $jk$-bigraded
$\Omega$-$\Omega$-bimodules with zero differential
$$\begin{array}{ll}
\Lambda^{1\diamond \bullet} & \cong \Omega \langle p \rangle[1-p] \oplus\Theta^\sigma \\
& \cong \Omega \langle p \rangle[1-p] \oplus(\Theta^\sigma)^\gamma \\
& \cong \mathbb{H}(X) \\
& \cong \mathbb{H}(\Psi).
\end{array}$$
\end{enumerate}
\end{lem}

Before proving the lemma, let us pause a moment on the statement of
part \eqref{26iv}. Implicit is the existence of a $\Omega$-$\Omega$-bimodule
structure on $\mathbb{H}(X)$. The reason for the existence of such a
structure is the action of $\mathbb{H} \End_\Omega(X)$ on
$\mathbb{H}(X)$, and the isomorphism $\mathbb{H} \End_{\Omega}(X)
\cong \Omega$ of part \eqref{26ii}.

\proof \eqref{26i} and \eqref{26ii}. As the image of the autoequivalence $\bft \otimes_{\bfc} -$ of $D^b(\bfc)$  under Koszul duality (and grading reversal), the functor $\Psi \otimes_{\Omega} -$ induces an autoequivalence of $D^b(\Omega)$. The image of the regular module
${}_\Omega\Omega$ under a derived auto equivalence
$\circlearrowright D^b(\Omega)$ is necessarily a tilting complex
whose endomorphism ring is isomorphic to $\Omega$.

\eqref{26iii} The isomorphism as $j$-graded left $\Omega$-modules is immediate from the definition of $X$ and Lemma \ref{ses}. The $k$-grading is forced upon us by the requirement that the differential has $k$-degree one, the existing $k$-grading on $\Omega$ and the fact that the term $\Omega e_{p-l}$ in homological position zero does indeed appear without any shifts as the image of the (unshifted) simple $\bfc$-module $L(p-l)$ under Koszul duality. Indeed, since $x^{p-l}$ has $k$-degree $p-l$, but the map given by right multiplication has to have $k$-degree one, the top of $\Omega e_p$ in the summand $\Omega e_p \to \Omega e_l$ of $X$ has to be concentrated in $k$-degree $p-l-1$, meaning that (since our $k$-grading convention follows that of the derived category, see Section \ref{gradings}) the corresponding $k$-graded version of the summand is $\Omega e_p[1+l-p] \to \Omega e_l$. The top of the homology given by the kernel of these maps (in every summand) is concentrated in $k$-degree $p-1$, hence the component $\Omega$ of $\mathbb{H}(X)$ is shifted by $[1-p]$.

\eqref{26iv} The first isomorphism follows from definition of $\Lambda$. Indeed $\Lambda^{1\diamond \bullet} = \Omega \otimes \zeta \oplus\Theta^\sigma \otimes 1$ and the isomorphism follows from the fact that $\zeta$ is a formal variable placed in $jk$-degree $(p,p-1)$ (noting that the degree of the top of $\Omega \langle p \rangle[1-p]$ is concentrated in $j$-degree $p$ and $k$-degree $p-1$, see Section \ref{gradings}).
The second isomorphism holds by Lemma \ref{inner}.

The remaining isomorphisms as left $jk$-bigraded $\Omega$-modules follow from (iii) and Lemma
\ref{specitilt}. To establish these isomorphisms of
$\Omega$-$\Omega$-bimodules, we need to examine the right action of $\Omega$ on $X$.

There are two obvious endomorphisms $\tilde y$ and $\tilde x$ of $X$ given by the commutative diagrams

\begin{equation}\label{y1}\xymatrix{
\Omega e_p \langle p-l \rangle \ar^{\cdot x^{p-l}}[rr] \ar_1[d] && \Omega e_{l} \ar^{\cdot x}[d] \\
\Omega e_p \langle p-l+1 \rangle \ar^{\cdot x^{p-l+1}}[rr] && \Omega e_{l-1} }\end{equation}
\begin{equation}\label{y2}\xymatrix{
\Omega e_{p} \langle p-1 \rangle \ar^{\cdot x^{p-1}}[rr] \ar[d]_1 && \Omega e_{1} \ar[d] \\
\Omega e_p \langle p \rangle \ar[rr] && 0 }\end{equation}

and

\begin{equation}\label{x1}\xymatrix{
\Omega e_p \langle p-l \rangle \ar^{\cdot x^{p-l}}[rr] \ar_{\cdot xy}[d] && \Omega e_{l} \ar^{\cdot y}[d] \\
\Omega e_p \langle p-l-1 \rangle \ar^{\cdot x^{p-l-1}}[rr] && \Omega e_{l+1} }\end{equation}
\begin{equation}\label{x2}\xymatrix{
\Omega e_p \langle p \rangle \ar[rr] \ar_{\cdot xy}[d]&& 0 \ar[d] \\
\Omega e_p \langle p-1 \rangle \ar^{\cdot x^{p-1}}[rr] && \Omega e_{1} }\end{equation}

respectively. On the kernel of the differential which is as a left module isomorphic to $\Omega$, it is easy to see that these endomorphisms satisfy the relations of $\Omega$. Indeed, in \eqref{y1} the kernel of the first row is $\Omega e_{p-l}$ and the kernel of the second row is $\Omega e_{p-l+1}$ and the induced action on there is right multiplication by $y$. Similary, the induced action on homology of $\tilde x$ is multiplication by $x$ on the kernel of the differential. 

To obtain our twist by $\gamma$, i.e.\ our claim on signs of this action, we phrase this purely in terms of dg bimodules.
As a dg-bimodule, the first row in \eqref{y1} is equal to $\Omega e_p \langle p-l \rangle [1+l-p]\oplus \Omega e_{l}$ with differential given by
$$\xymatrix{\Omega e_p \langle p-l \rangle [1+l-p]\oplus \Omega e_{l} \ar^{{\footnotesize \left(\begin{array}{cc}0&0\\ \cdot x^{p-l}[p-l-1]&0\end{array}\right)}}[rrrr]&&&&\Omega e_p \langle p-l \rangle [1+l-p]\oplus \Omega e_{l} }$$
where the shifts in $k$-grading are again forced on us taking into account that $x$ has $k$-degree $1$ and the differential needs to have $k$-degree $1$.
So \eqref{y1} translates to 
$$\xymatrix{\Omega e_p \langle p-l \rangle [1+l-p]\oplus \Omega e_{l}\ar_{ {\footnotesize \left(\begin{array}{cc}[-1]&0\\ 0&0\end{array}\right)}}[dd]
\ar^{{\footnotesize \left(\begin{array}{cc}0&0\\ \cdot x^{p-l}[p-l-1]&0\end{array}\right)}}[rrrr]&&&&\Omega e_p \langle p-l \rangle [1+l-p]\oplus \Omega e_{l}\ar^{
{\footnotesize \left(\begin{array}{cc}0&0\\ 0&\cdot x \end{array}\right)}
}[dd] \\&&&& \\
\Omega e_p \langle p-l+1\rangle [l-p]\oplus \Omega e_{l} \ar^{{\footnotesize \left(\begin{array}{cc}0&0\\ \cdot x^{p-l+1}[p-l]&0\end{array}\right)}}[rrrr]&&&&\Omega e_p \langle p-l+1\rangle [l-p]\oplus \Omega e_{l} }$$

In order to make the dg-formalism $$(a.m)d = a.(m)d, \quad (m.a)d = (-1)^{|a|_k} m(d).a$$
described in Section \ref{gradings} work, we see that for our generator for our generator $\tilde y  \cong \End_{D^b(\Omega)}(X)$, we need (setting $m=e_p$ and $a=\tilde y$), that \begin{equation}\label{sign1}e_p\tilde y \cdot x^{p-l}[p-l-1] =  - e_p \cdot x^{p-l+1}[p-l] \cdot \tilde y,\end{equation}  and similarly for the action of $\tilde x$. Therefore, choosing the isomorphism  $\End_{D^b(\Omega)}(X) \to \Omega$ to send $\tilde x$ and $\tilde y$ to $x$ and $y$ respectively, we need to twist the action on $\Theta^{\sigma}$ by $\gamma$ on the right, i.e.\ $\tilde y$ will correspond to right multiplication by $-x$ on $\Theta$, which then, picking up the sign from \eqref{sign1}, will make our diagrams commute.
\endproof

\subsection{The homology $\mathbb{H}(\Psi^{\otimes_\Omega i})$.}
In order to prove that $\Lambda$ is isomorphic to
$\mathbb{H}(\bbT_{\Omega}(\Psi))$, we now need to consider higher $i$-degrees, in other words we need to compute the homology of $\Psi^{\otimes_\Omega i}$.

Consider the complex $X_i$ of $\Omega$-modules, obtained by splicing
together $n$ copies of $X$. To be more precise, let $X_i$ denote the
complex
$$X^1 \langle (i-1)p \rangle \rightarrow W^1 \langle (i-2)p \rangle \rightarrow ... \rightarrow W^1 \langle p \rangle \rightarrow W^1 \rightarrow W^0$$
whose differentials are obtained via the composition
$$\xymatrix{
W^1 \langle p \rangle \ar[rr] \ar[dr]_{d_W \langle p \rangle} & & W^1  \\
& W^0 \langle p \rangle \ar[ur]_{\delta}.&\\}$$ 
In other words, $X_i$ is the direct sum of 
$$\Omega e_p \langle ip \rangle \rightarrow 0  \rightarrow ... \rightarrow 0$$
where the nonzero term is in the homological position $-i$, and
$$\bigoplus_{l=1}^{p-1}(\Omega e_p \langle l+ (i-1)p \rangle \rightarrow \Omega e_p \langle  p-l+(i-2)p \rangle \rightarrow ... \rightarrow \Omega e_p\langle (p-l) +p \rangle \rightarrow \Omega e_p\langle l \rangle  \rightarrow \Omega e_{p-l})$$
if $i$ is odd or 
$$\bigoplus_{l=1}^{p-1}(\Omega e_p \langle l+ (i-1)p \rangle \rightarrow \Omega e_p \langle p- l+(i-2)p \rangle \rightarrow ... \rightarrow \Omega e_p\langle l +p \rangle \rightarrow \Omega e_p\langle p-l \rangle  \rightarrow \Omega e_{l})$$
if $i$ is even.
Here again the maps right multiplication by $x^l$ respectively $x^{p-l}$ for the rightmost arrow and right multiplication by 
the appropriate power of $xy$ dictated by the grading for the remaining arrows.

We define $\Psi_i$ to
be $\Psi^{\otimes_{\Omega} i}$.

\begin{lem} \label{Lambda2}
\begin{enumerate}[(i)]
 \item\label{27i} We have a quasi-isomorphism of differential bigraded
$\Omega$-modules
$$X_i \cong \Psi_i.$$
\item \label{27ii} We have isomorphisms of $jk$-bigraded $\Omega$-$\Omega$-bimodules
with zero differential
$$\begin{array}{ll}
\Lambda^{i \diamond \bullet} &
\cong \Omega \langle ip \rangle[i(1-p)] \oplus \Theta^{\sigma} \langle (i-1)p \rangle[(i-1)(1-p)]
\oplus ...\oplus \Theta^{\sigma^i} \langle 0 \rangle[0] \\
& \cong \Omega \langle ip \rangle[i(1-p)] \oplus \Theta^{\sigma\gamma} \langle (i-1)p \rangle[(i-1)(1-p)]
\oplus ...\oplus \Theta^{(\sigma\gamma)^i} \langle 0 \rangle[0] \\
& \cong \mathbb{H}(X_i) \\
& \cong \mathbb{H} (\Psi_i).
\end{array}$$
\item \label{27iii} For a suitable choice of such isomorphisms, the multiplication
map
$$\Lambda^i \otimes \Lambda^{i'} \rightarrow \Lambda^{i+i'}$$
corresponds to $\mathbb{H}$ applied to the multiplication map
$$\Psi_i \otimes \Psi_{i'} \rightarrow \Psi_{i+i'}.$$
\end{enumerate}

\end{lem}

\proof The of the functor $\Psi \otimes_{\Omega} -$ (as an autoequivalence of the bounded derived category of $\Omega$-modules) is to shift the projective $\Omega e_{p}$ by $1$, whilst
taking $\Omega e_{l}$ to a complex $\Omega e_{p} \rightarrow \Omega
e_{l}$ concentrated in degrees $-1$ and $0$ whose differential has image covering the image of all
possible maps $\Omega e_{p} \rightarrow \Omega e_{l}$. If we iterate
this construction $i$ times in a $jk$-graded setting, we obtain precisely
$X_i$. This establishes \eqref{27i}.

\eqref{27ii} The first isomorphism comes from the definition of $\Lambda$.
Indeed, we have $$\Lambda^{i \diamond \bullet}  = \Omega \otimes \zeta^i 
 \oplus \Theta^{\sigma}\otimes \zeta^{i-1}\oplus ...\oplus \Theta^{\sigma^i} \otimes 1$$ and again recalling that the formal variable $\zeta$ is placed in $jk$-degree $(p,p-1)$ and that by our conventions in Section \ref{gradings} these are indeed the degrees in which the components of homology are concentrated, we obatin the claim. The second isomorphism again follows from Lemma \ref{inner}.

For the third isomorphism, consider the homology
concentrated in the middle term of

\xymatrix{
W^1 \langle (l+1)p \rangle \ar[rr] \ar[dr]_{d_W \langle (l+1)p \rangle} & & W^1 \langle lp \rangle \ar[rr]\ar[dr]_{d_W \langle lp \rangle} & & W^1\langle (l-1)p \rangle \\
& W^0 \langle (l+1)p \rangle \ar[ur]_{\delta\langle lp \rangle } &&
W^0 \langle lp \rangle \ar[ur]_{\delta\langle (l-1)p \rangle} &\\}

which is given by $$Q = \frac{\Ker d_W \langle lp \rangle}{\im
(\delta\langle lp \rangle \circ d_W \langle (l+1)p \rangle )}. $$
Note that the subsequence 
$$ W^0 \langle (l+1)p \rangle \rightarrow
W^1 \langle lp \rangle \rightarrow W^0 \langle lp \rangle $$ 
which has summands of the form
$$\Omega e_l \langle (l+1)p \rangle \rightarrow \Omega e_p \langle lp+ l \rangle \rightarrow \Omega e_{p-l} \langle lp \rangle  $$
is exact on the left and in the middle by Lemma \ref{ses}, so in fact
$$Q \cong \Coker d_W \langle (l+1)p \rangle. $$ This was computed in
Lemma \ref{ses} as a left $\Omega$-module. 

The claim on the $k$-grading follows in the same way as in Lemma  \ref{Lambda1}\eqref{26iv}, whose computations we also follow
to establish a bimodule structure.
Again, we know that we have an isomorphism $\End_{D^b(\Omega)}(X_i)\cong \Omega$, and we have obvious generators (in the example where $i$ is odd and $1 \leq l \leq p-2$)

\begin{equation}\label{y1i}\xymatrix@C=10pt{
\Omega e_p \langle l+ (i-1)p \rangle \ar[r]\ar^{\cdot 1}[d]& \Omega e_p \langle  p-l+(i-2)p \ar^{\cdot xy}[d] \rangle & \cdots &  \Omega e_p\langle l \rangle \ar^{\cdot 1}[d] \ar[r]&\Omega e_{p-l}\ar^{\cdot x}[d]\\
\Omega e_p \langle l+1+ (i-1)p \rangle \ar[r]& \Omega e_p \langle  p-l-1+(i-2)p \rangle & \cdots &  \Omega e_p\langle l+1 \rangle  \ar[r]&\Omega e_{p-l-1}
}
\end{equation}

and 

\begin{equation}\label{x1i}\xymatrix@C=10pt{
\Omega e_p \langle l+ (i-1)p \rangle \ar[r]\ar^{\cdot xy}[d]& \Omega e_p \langle  p-l+(i-2)p \ar^{\cdot 1}[d] \rangle & \cdots &  \Omega e_p\langle l \rangle \ar^{\cdot xy}[d] \ar[r]&\Omega e_{p-l}\ar^{\cdot y}[d]\\
\Omega e_p \langle l-1+ (i-1)p \rangle \ar[r]& \Omega e_p \langle  p-l+1+(i-2)p \rangle & \cdots &  \Omega e_p\langle l-1 \rangle  \ar[r]&\Omega e_{p-l+1}
}
\end{equation}

Again choosing the identification with $\Omega$ in such a way that it coincides with natural right multiplication on the component in homology isomorphic to $\Omega$ appearing in the left-most homological position and obeying the sign rule we pick up a twist by $\sigma\gamma$ in every step.

The last isomorphism follows from \eqref{27i}.

\eqref{27iii} Here we write $\Lambda^i = \oplus_{j,k \in \mathbb{Z}}
\Lambda^{ijk}$. By its definition as a tensor algebra, the algebra $\Lambda$ is quadratic with respect to
the $i$-grading, with generators $\Lambda^0  = \Omega$ in degree
$0$, generators $\Lambda^1 = \Omega\zeta \oplus \Theta^\sigma$ in degree
$1$, and relations
$$R \subset (\zeta\Omega \otimes_\Omega \Theta^\sigma) \oplus (\Theta^\sigma \otimes_\Omega \zeta\Omega) \subset \Lambda^1 \otimes_{\Lambda^0} \Lambda^1$$
given by the image of the composition
$$\xymatrix{
\Theta^\sigma \ar[rr]^{a \mapsto (a,-a)} & & \Theta^\sigma \oplus
\Theta^\sigma \ar[rrr]^{(a,b) \mapsto (\zeta \otimes a, b \otimes \zeta)} &
& & (\zeta\Omega \otimes_\Omega \Theta^\sigma) \oplus (\Theta^\sigma
\otimes_\Omega \zeta\Omega) },$$ 
which simply says that the formal variable $\zeta$ and the bimodule $ \Theta^\sigma$ commute.

It follows that to prove \eqref{27iii}, it is sufficient to check that the hexagon
$$\xymatrix{
& \ar[dl]_\sim \Theta^\sigma \ar[dr]^\sim & \\
\Theta^\sigma \otimes_\Omega \Omega \ar[d]_\wr & & \Omega \otimes_\Omega \Theta^\sigma \ar[d]^\wr \\
\mathsf{H}^0(\Psi_1) \otimes_\Omega \mathsf{H}^1(\Psi_1) \ar[dr]_\mu & & \ar[dl]^\mu \mathsf{H}^1(\Psi_1) \otimes_\Omega \mathsf{H}^0(\Psi_1) \\
& \mathsf{H}^1(\Psi_2) & }$$ commutes, where $\mu$ denotes the
multiplication map in $\mathbb{H}\mathbb{T}_\Omega(\Psi)$ and $\mathsf{H}^q$ denotes the homology in the $q$th position of $\Psi_i$ regarded as a complex of $\Omega$-$\Omega$-bimodules in $D^b(\Omega)$ (which is not the same as taking the $q$th homology $H^q$ with respect to the $k$-grading (cf.\ Section \ref{gradings}) as again $\Omega$ is not concentrated in $k$-degree zero, so we have two different decompositions $\mathbb{H} = \bigoplus_q \mathsf{H}^q$ and $\mathbb{H} = \bigoplus_q H^q$).
Hence we have bimodule isomorphisms 
\begin{equation}\label{psiisos}\begin{split}
\mathbb{H}(\Psi) &=  \mathsf{H}^1(\Psi) \oplus  \mathsf{H}^0(\Psi)\cong \Omega \langle p \rangle [1-p]\oplus  \Theta^{\sigma\gamma} \\
\mathbb{H}(\Psi_2)& = \mathsf{H}^2(\Psi_2) \oplus \mathsf{H}^1(\Psi_2) \oplus \mathsf{H}^0(\Psi_2)\cong \Omega \langle 2p \rangle [2-2p] \oplus  \Theta^{\sigma\gamma} \langle p \rangle [1-p] \oplus \Theta.\end{split}\end{equation}
Our explicit description of the action of $\Omega$ on $X_i$ also gives fixed bimodule isomorphisms

\begin{equation}\label{xisos}\begin{split}
\mathbb{H}(X) &=  \mathsf{H}^1(X) \oplus  \mathsf{H}^0(X)\cong \Omega\langle p \rangle [1-p]\oplus  \Theta^{\sigma\gamma} \\
\mathbb{H}(X_2) &= \mathsf{H}^2(X_2) \oplus \mathsf{H}^1(X_2) \oplus \mathsf{H}^0(X_2)\cong \Omega\langle 2p \rangle [2-2p] \oplus  \Theta^{\sigma\gamma}\langle p \rangle [1-p] \oplus \Theta.\end{split}\end{equation}

Let us now explain how once we have fixed our
quasi-isomorphism $\phi$ between $X$ and $\Psi$, we can define an associated fixed isomorphism
$$\mathbb{H}(X_2) \cong \mathbb{H}(\Psi_2)$$
hence fixing bimodule isomorphisms in \eqref{psiisos} as induced by those in \eqref{xisos}.

Indeed we have an isomorphism
$$\mathbb{H}(\Psi \otimes \phi): \mathbb{H}(\Psi \otimes_\Omega X) \rightarrow \mathbb{H}(\Psi_2);$$ and the
isomorphisms
$$\mathbb{H}(\Psi \otimes_\Omega X) = \mathbb{H} \left( \Psi \otimes_\Omega \left( (\Omega e_p \otimes_F e_p FQ e \rightarrow \Omega e) \oplus (\Omega e_p \langle p \rangle \rightarrow 0) \right) \right),$$
$$\mathbb{H}(\phi e_l): \mathbb{H}(Xe_l) \cong \mathbb{H}(\Psi \otimes_\Omega \Omega e_l)$$
imply an isomorphism
$$\mathbb{H}\left( (Xe_p \otimes_F e_p FQ e \rightarrow X e) \oplus (X e_p \langle p \rangle \rightarrow 0) \right) \cong \mathbb{H}(\Psi \otimes_\Omega X),$$
Our explicit description of $X$, of $X_2$, and of the right action of $\Omega$ on $X$,
induces an isomorphism
$$\mathbb{H}(X_2) \rightarrow \mathbb{H}\left( (Xe_p \otimes_F e_p FQ e \rightarrow X e) \oplus (X e_p \langle p \rangle \rightarrow 0) \right).$$
Composing these isomorphisms gives us our fixed isomorphism
$\mathbb{H}(X_2) \cong \mathbb{H}(\Psi_2)$.

To  establish our commuting hexagon we now proceed to examine
in detail the effect of right multiplication by $\mathbb{H}(\Psi)$ upon fixing the above identifications, and by compatibility of our morphisms we can do this on the level of $X$, with our very explicit description.
We have an explicit isomorphism $$\mathbb{H}(X_2) \cong \mathbb{H} \left( X \otimes_\Omega \left( (\Omega e_p \otimes_F e_p FQ e \rightarrow \Omega e) \oplus (\Omega e_p \langle p \rangle \rightarrow 0) \right) \right)$$
and in this formulation, given our analysis in Lemma \ref{Lambda1}\eqref{26iv}, we can see that an element $\omega\in\mathsf{H}^1(X) \cong \Omega\langle p \rangle [1-p]$ in the  second tensor factor (which is isomorphic to $\Omega$), will act on the right of an element $\theta \in \mathsf{H}^0(X)\cong\Theta^{\sigma\gamma}$ by the formula $\theta \cdot \omega = \theta \gamma\sigma(\omega)$. Similarly and element $\theta' \in \mathsf{H}^0(X)\cong\Theta^{\sigma\gamma}$ in the second tensor factor will act on the right on and element $\omega' \in \mathsf{H}^1(X) \cong \Omega\langle p \rangle [1-p]$ simply by the formula $\omega' \cdot \theta' = \omega\theta$.

Note that the  canonical bimodule isomorphism $\Theta^{\sigma \gamma} \to \Theta^{\sigma}$ implied by Lemma \ref {inner} is given by the map $\theta \mapsto \theta g$, where $g$ was the self-inverse element inducing the inner automorphism $\gamma$.

This implies that the hexagon of fixed bimodule isomorphisms (where we have used our identifications $\mathsf{H}^0(X)\otimes_\Omega\mathsf{H}^1(X)\cong \Theta^{\sigma\gamma} \otimes_\Omega \Omega $,  $\mathsf{H}^1(X)\otimes_\Omega\mathsf{H}^0(X)\cong \Omega \otimes_\Omega \Theta^{\sigma\gamma}$ and   $\mathsf{H}^1(X_2) \cong \Theta^{\sigma}$)

$$\xymatrix{
& \theta \in \ar@{|->}[dl] & \ar[dl]_\sim \Theta^{\sigma} \ar[dr]^\sim \ar[ddd]^{\cdot g} & \ni \theta  \ar@{|->}[dr]& \\
\theta \otimes 1	 \in  \ar@{|->}[d] & \Theta^{\sigma} \otimes_\Omega \Omega \ar[d]_\wr & &
\Omega  \otimes_\Omega \Theta^{\sigma} \ar[d]^\wr &  \ni1 \otimes \theta \ar@{|->}[d] \\
\theta g \otimes 1 \in \ar@{|->}[dr] &\mathsf{H}^0(X)\otimes_\Omega\mathsf{H}^1(X)\ar[dr] & & \ar[dl] \mathsf{H}^1(X)\otimes_\Omega\mathsf{H}^0(X)& \ni 1 \otimes \theta g \ar@{|->}[dl] \\
& \theta g \sigma\gamma(1) =\theta g\in & \mathsf{H}^1(X_2) & \ni \theta g }$$

commutes, which completes the proof of the Lemma.
\endproof

\subsection{Proof of Theorem \ref{Lambdacalculation}.}
The isomorphism
$\mathbb{H}( \mathbb{T}_\Omega(\Psi)) \cong \Lambda$ as an $\Omega$-$\Omega$-bimodule comes from Lemma
\ref{Lambda2}\eqref{27ii}. It is an algebra isomorphism, thanks to Lemma
\ref{Lambda2}\eqref{27iii}.
Recall the gradings on both algebras:

The $i$-grading on the tensor algebra $\bbT_{\Omega}(\Psi)$ has
$\Omega$ in $i$-degree $0$ and $\Psi$ in $i$-degree $1$, while the $j$-grading is inherited by regarding $\Psi$ as a differential
$jk$-bigraded $\Omega$-module, where the $j$- and $k$-gradings on $\Omega$ are both by path length.

We have $\Lambda = \bbT_{\Omega}(\Theta^{\sigma}) \otimes F[\zeta]$;
we let $\Omega$ have $i$-degree $0$ and let $\zeta$ and $\Theta^\sigma$
have $i$-degree $1$; the $j$-grading on $\Lambda$ is obtained by grading $\Omega$
and $\Theta^{\sigma}$ by path length and placing $\zeta$ in degree $p$;
finally the $k$-grading is obtained by grading $\Omega$ and $\Theta^{\sigma}$ by path length
and placing $\zeta$ in degree $p-1$.

These gradings match up under the isomorphism, hence $\Lambda \cong \mathbb{H}
\bbT_{\Omega}(\Psi)$ as $ijk$-trigraded algebras.\qed

\section{Computing $\bfy$.} \label{computation}

Here we apply the theory of the preceding sections to compute $\bfy$.

\emph{Proof of Theorem \ref{Yoneda}}. We have algebra isomorphisms

$$\begin{array}{lll}
\bfy_q  & \cong \mathbb{H} \bbO_{F,0} \bbE \bbO_{\bfc,\bft}^q(F,F) & \mbox{Theorem \ref{GL2theorem}} \\
 & \cong \mathbb{H} \bbO_{F,0} \bbE \bbO_{\bfc,\tilde\bft}^q(F,F) & \mbox{Lemma \ref{cqiso}} \\
& \cong \mathbb{H} \bbO_{F,0} \bbO_{\bfc^!,\bft^!}^q \bbE(F,F) & \mbox{Corollary \ref{quasiiterated}} \\
& \cong \mathbb{H} \bbO_{F,0} \bbO_{\bfc^!,\bft^!}^q(F,F) & \mbox{$\bbE(F,F)=(F,F)$} \\
& \cong \mathbb{H} \bbO_{F,0} \bbO_{\Omega, \Psi}^q(F,F) & \mbox{Corollary \ref{minustopluscor}}\\
& \cong \mathbb{H} \fO_F \fO_{\bbT_{\Omega}(\Psi)}^q(F[z]) & \mbox{Corollary \ref{compare}} \\
& \cong \fO_{\mathbb{H} F} \fO_{\mathbb{H}\bbT_{\Omega}(\Psi)}^q \mathbb{H}(F[z]) & \mbox{Lemma \ref{commute}} \\
& \cong \fO_F \fO_\Lambda^q(F[z]) & \mbox{Theorem \ref{Lambdacalculation}} \\
& \cong \lambda_q & \mbox{definition of $\lambda_q$} \end{array}$$

Sending $q$ to $\infty$, we obtain $\bfy \cong \lambda$ as required.
\endproof

\normalfont

{\sc Vanessa Miemietz, Will Turner}\\
School of Mathematics, University of East Anglia, Norwich, NR4 7TJ, UK, \\{\tt v.miemietz@uea.ac.uk}\\
Department of Mathematics, University of Aberdeen, Fraser Noble Building, King's College, Aberdeen AB24 3UE, UK, {\tt w.turner@abdn.ac.uk}.

\begin{thebibliography}{00}
\normalfont

\bibitem{BGS} A. Beilinson, V. Ginzburg, W. Soergel, \emph{Koszul duality patterns in representation theory}, J. Amer. Math. Soc.  9  (1996),  no. 2, 473--527.

\bibitem{CPS}  E. Cline, B. Parshall and L. Scott, \emph{Finite-dimensional algebras and highest weight categories.}
J. Reine Angew. Math. 391 (1988), 85--99.


\bibitem{EH} K. Erdmann, A. Henke, \emph{On Ringel duality for Schur algebras},
Math. Proc. Cambridge Philos. Soc. 132 (2002), no. 1, 97--116. 

\bibitem{Green} J. A. Green, \emph{Polynomial representations of ${\rm
GL}_{n}$.} Lecture Notes in Mathematics, 830. Springer, Berlin,
1980.

\bibitem{Ke}
B. Keller, \emph{On differential graded categories}, International Congress of
  Mathematicians. Vol. II, Eur. Math. Soc., Z\"urich, 2006, pp.~151--190.

\bibitem{MT1} V. Miemietz, W. Turner, \emph{Rational representations of $GL_2$}, Glasgow J. Math. 53 (2011), no.2, 257--275.

\bibitem{MT2} V. Miemietz, W. Turner, \emph{Homotopy, Homology and $GL_2$}, Proc. London Math. Soc. doi:10.1112/plms/pdp040.

\bibitem{AP} A. Parker, \emph{Higher extensions between modules for $SL_2$},  Adv. Math.  209  (2007),  no. 1, 381--405.

\bibitem{Rickard}  J. Rickard, \emph{Derived equivalences as derived functors.} J. London Math. Soc. (2) 43 (1991), no. 1, 37-48. 

\bibitem{RR} R. Rouquier, \emph{Derived equivalences and finite
dimensional algebras}, Proceedings of the International Congress of
Mathematicians (Madrid, 2006 ), vol II, pp. 191-221, EMS Publishing
House, 2006.





\end{thebibliography}
\end{document}